\documentclass[12pt]{article}

\input{xypic.tex} \xyoption{all}
\usepackage{amsthm,amssymb,latexsym}
\usepackage{graphics}

\newcommand{\seta}{\rightarrow} 
\newcommand{\real}{\mathbb{R}} \newcommand{\cpx}{\mathbb{C}} 
\newcommand{\zed}{\mathbb{Z}} 
\newcommand{\p}{\mathbb{P}}
\newcommand{\st}{\mathbb{T}} \newcommand{\bl}{\mathbb{L}}
\newcommand{\f}{\mathbb{F}} \newcommand{\m}{\mathbb{M}}
\newcommand{\mi}{\mathbb{M}^{\rm I}} 
\newcommand{\cp}{\mathfrak{P}} \newcommand{\cm}{\mathfrak{M}}
\newcommand{\cf}{\mathfrak{F}} 
 
\newcommand{\sm}{\mathsf{m}} \newcommand{\sn}{\mathsf{n}}
 
\newcommand{\id}{\mathbf{1}}

\newtheorem{theorem}{Theorem}
\newtheorem{proposition}[theorem]{Proposition}
\newtheorem{lemma}[theorem]{Lemma}
\newtheorem*{definition}{{\bf Definition}}

\begin{document}

\title{Quantum instantons \\ with classical moduli spaces}
\author{Igor Frenkel \\ Yale University \\ Department of Mathematics \\
10 Hillhouse Avenue \\ New Haven, CT 06520-8283 USA\\ \\
Marcos Jardim \\ University of Massachusetts at Amherst \\
Department of Mathematics and Statistics \\
Amherst, MA 01003-9305 USA }

\maketitle

\begin{abstract}
We introduce a quantum Minkowski space-time based on the quantum group
$SU(2)_q$ extended by a degree operator and formulate a quantum version
of the anti-self-dual Yang-Mills equation. We construct solutions of
the {\em quantum} equations using the {\em classical} ADHM linear data,
and conjecture that, up to gauge transformations, our construction
yields all the solutions. We also find a deformation of Penrose's
twistor diagram, giving a correspondence between the quantum Minkowski
space-time and the classical projective space $\p^3$.
\end{abstract}

\newpage
\baselineskip18pt

%---------------------------------------------------------------
%---------------------------------------------------------------

\section*{Introduction} \label{intro}

Our present view of the mathematical structure of 
space-time was first formulated by H. Minkowski
in \cite{Mi}, based on Einstein's discovery of 
Special Relativity. Since then, many mathematicians
and physicists have tried to develop various 
generalizations of
Minkowski's concept of space-time, 
in different directions.

One such generalization was
strongly advocated by R. Penrose \cite{P}. He studied
the conformal compactification of the complexified
Minkowski space-time, denoted by $\m$, and the
associated space of null straight lines, denoted by
$\p$, showing that various equations from Mathematical
Physics (e.g. Maxwell and Dirac equations) over 
$\m$ can be transformed into natural holomorphic 
objects over $\p$. 

One of the most celebrated examples of the Penrose 
programme was the solution of the anti-self-dual 
Yang-Mills (ASDYM) equations by M. Atiyah, V. Drinfeld, 
N. Hitchin and Yu. Manin. These authors described explicitly 
the moduli space ${\cal M}^{\rm reg}(n,c)$ of finite action
solutions of the ASDYM equation over $\real^4$, usually
called {\em instantons}, of fixed rank $n$ and charge $c$
in terms of some linear data \cite{ADHM}. In particular,
they constructed an instanton associated to any linear data,
and proved, conversely, that every instanton can be obtained 
in this way, up to gauge transformations.

The Penrose approach has also been successfully applied to 
massless linear equations, the full Yang-Mills equation, 
self-dual Einstein equation, etc.
(see \cite{WW} and the references therein).

Another generalization of the Minkowski space-time based on
noncommutative geometry was proposed by A. Connes \cite{C}.
The first step in this approach is the replacement of
geometric objects by their algebras of functions and the reformulation
of various geometric concepts in algebraic language. The next
step is the deformation of the algebraic structures and the
introduction of noncommutativity. The {\it algebraic} structures so
obtained are no longer associated with the original geometry,
being regarded as {\it geometric} structures on a ``noncommutative space''. 

The simplest example of noncommutative space-time can be obtained
by replacing the commutative algebra $\real^4$ by the Heisenberg algebra.
Recently, N. Nekrasov and A. Schwarz \cite{NS} have defined instantons
on such noncommutative $\real^4$ and argued that the corresponding moduli
space can be parametrized by a modification of the Atiyah, Drinfeld,
Hitchin and Manin (ADHM) linear data.

Another source for the deformation of the Minkowski space-time
emerged with the discovery of {\em quantum groups} by V. Drinfeld
\cite{Dr} and M. Jimbo \cite{J}. It quickly led to the notion of
a quantum Minkowski space-time and various related structures;
some of the early papers are \cite{CSSW,Mj,SWZ}. In particular,
a quantum
version of the ASDYM equation was studied in \cite{Z}. 

One expects that the linear data of Atiyah, Drinfeld, Hitchin and Manin
should also be deformed in order to yield solutions of the quantum
ASDYM equation. Thus the quantum Minkowski space-time and related
structures, though mathematically sensible, seem rather dubious in
terms of physical applications. In fact, the quantum deformation
destroys the classical symmetry groups and the possible reconstruction
of these symmetries is far from apparent.

In the present paper we observe a new phenomenon, which goes against
the typical intuition related to quantum deformations of the Minkowski
space-time and various equations on it. We show that there exists a
natural quantum deformation of the Minkowski space-time (in fact, of 
the whole compactified complexified Minkowski space $\m$, along with
its real structures) and the ASDYM equation, such that {\em the moduli
space of quantum instantons is naturally parameterized by the classical,
non-deformed ADHM data}.

Starting from the classical ADHM data, we will explicitly construct
solutions of the quantum ASDYM equation. Furthermore, we conjecture that
our construction yields all the solutions up to gauge equivalence. We
hope to prove this conjecture using a generalization of the Penrose
twistor transform, which takes the {\em quantum} Minkowski space-time
to the {\em classical} space of straight lines. Such procedure thus
realizes Penrose's dream, who regarded light rays as more fundamental
than points in space-time: in our construction, space-time is being
deformed while the space of light rays is kept fixed.

Moreover, this phenomenon of quantum equations with classical
solutions does not seem to be restricted only to the
ASDYM equation, but it is also steadfast for massless linear equations,
full Yang-Mills equations, etc. This opens new venues for physical
applications of the quantum Minkowski space-time here proposed: the
classical symmetry groups can be restored if one considers {\em all}
quantum deformations of the classical Minkowski space-time, since the
spaces of solutions of the quantum equations admit natural identifications.

Our constructions are based on the theory of the quantum group $SU(2)_q$
extended by a natural degree operator, which makes possible for
a surprising construction of solutions of the quantum ASDYM equation from
the classical ADHM data. The extended quantum group $\widetilde{SU}(2)_q$
has functional dimension $4$, the same as the quantum space-time,
and we show that the relation between them extends to deep structural levels.

To formulate the ASDYM equation we need a theory of exterior
forms on quantum Minkowski space-time. This is derived from the differential
calculus on $SU(2)_q$ first developed in \cite{WZ,W}, with the addition of the
differential of the degree operator. The $R$-matrix formulation of the quantum
group $SU(2)_q$ \cite{FRT} and its exterior algebra \cite{S} allows us to
present the construction of the quantum connections even more compactly than in
the classical case, and efficiently verify the quantum ASDYM equation.

Our results bring us to the conclusion that {\em the correct notion for
a quantum Minkowski space-time is precisely the extended quantum group
$\widetilde{SU}(2)_q$}. This relation, which might seem artificial at the
classical level ($q=1$), is imposed on us by the mathematical structure
itself. We believe that future research will reveal the full potential of
this new incarnation of the Minkowski space-time.

\tableofcontents

%---------------------------------------------------------------
%---------------------------------------------------------------

\section{Quantum Minkowski space-time from $SU(2)_q$} \label{mink}

\subsection{Algebraic structures on Minkowski space-time}

We begin with some well known facts regarding Penrose's approach to the
Minkowski space-time for the convenience of the reader; further details
can be found in \cite{Ma,WW}.

Let $\st$ be a 4 dimensional complex vector space and consider 
$\m = \mathbf{G}_2(\st)$, the Grassmannian of planes in $\st$.
As usual in the literature, we will often refer to $\m$ as the
{\em compactified complexified Minkowski space}, since $\m$ can be
obtained via a conformal compactification of $M^4\otimes\cpx$,
where $M^4$ denotes the usual Minkowski space. 

The Grassmannian $\m$ can be realized as a quadric in $\p^5$ via the
{\em Pl\"ucker embedding}. More precisely, note that $\mathbf{P}(\Lambda^2\st)\simeq\p^5$,
and take homogeneous coordinates $[z_{rs}]$ for $r,s=1,2,3,4$ (where $z_{rs}=-z_{sr}$ is
the coefficient of $dz_r\wedge dz_s$). Then $\m$ becomes the subvariety of 
$\mathbf{P}(\Lambda^2\st)$  given by the quadric:
\begin{equation} \label{quadric}
z_{12}z_{34} - z_{13}z_{24} + z_{14}z_{23} = 0
\end{equation}

Let us now fix a direct sum decomposition of $\st$ into two
2-dimensional subspaces:
\begin{equation} \label{decomp}
\st = \bl \oplus \bl'
\end{equation}
Such choice induces a decomposition of the second exterior power as follows:
$$ \Lambda^2\st = \Lambda^2\bl \oplus \Lambda^2\bl' \oplus \bl\wedge\bl' $$
Now fix basis $\{e_1,e_2\}$ and $\{e_{1'},e_{2'}\}$ in $\bl$ and $\bl'$,
respectively. To help us keep track of the choices made, we will use indexes
$\{1,2,1',2'\}$ instead of $\{1,2,3,4\}$. Since the variables $z_{12}$ and
$z_{1'2'}$ will play a very special role in our discussion, we will
introduce the notation $D=z_{12}$ and $D'=z_{34}=z_{1'2'}$.
The quadric (\ref{quadric}) is then rewritten in the following way:
\begin{equation} \label{quad}
z_{11'}z_{22'}-z_{12'}z_{21'} = DD' 
\end{equation}

The decomposition (\ref{decomp}) also induces the choice of a {\em point at
infinity} in the compactified complexified Minkowski space $\m$. Let 
${\cal S}(x)$ denote the plane in $\st$ corresponding to the point
$x\in\m$ and $\ell$ denote the point in $\m$ corresponding to the plane
$\bl$. Consider the sets:
\begin{eqnarray}
\m^{\rm I} = \left\{ x\in\m\ |\ {\cal S}(x)\cap\bl=\{0\} \right\} & \ \ \ & 
{\rm complexified\ Minkowski\ space} \nonumber\\
C(\ell) = \left\{ x\in\m\ |\ \dim\left({\cal S}(x)\cap\bl\right)=1 \right\} & \ \ \ & 
{\rm light\ cone\ at\ infinity} \nonumber
\end{eqnarray}
Then clearly $\m=\m^{\rm I} \cup C(\ell) \cup \{\ell\}$, and $\m^{\rm I}$ is an affine
space, being isomorphic to $\cpx^4=M^4\otimes\cpx$. Moreover, we note
that the light cone at infinity $C(\ell)$ has complex codimension one in $\m$.

We will denote the local coordinates on $\m^{\rm I}$ 
by $x_{rs'}=z_{rs'}/D$, where $r,s=1,2$. They are related to the
Euclidean coordinates $x^k$ on $\m^{\rm I}$ in the following way:
\begin{equation} \label{coords1} \begin{array}{ccc}
x_{11'}  = x^1 - ix^4 & \ \ \ &
x_{12'} = -ix^2 - x^3 \\
x_{21'} = -ix^2 + x^3 & \ \ \ &
x_{22'} = x^1 + ix^4
\end{array}\end{equation}
We will denote by $E^4$ the real Euclidean space spanned by $x^k$.

Similarly, let $\ell'$ be the point in $\m$ corresponding to the plane $\bl'$.
It can be regarded as the {\em origin} in $\cpx^4$. We define:
$$ \m^{\rm J} = \left\{ x\in\m\ |\ {\cal S}(x)\cap\bl'=\{0\} \right\} $$
$$ C(\ell') = \left\{ x\in\m\ |\ \dim\left({\cal S}(x)\cap\bl'\right)=1 \right\} $$
so that $\m^{\rm J}$ is also a 4-dimensional affine space and  
$\m=\m^{\rm J} \cup C(\ell') \cup \{\ell'\}$. We will denote the 
local coordinates on $\m^{\rm J}$ by $y_{rs'}=z_{rs'}/D'$, where $r,s=1,2$.

It is important to note that even though the affine spaces $\m^{\rm I}$ and
$\m^{\rm J}$ do not cover the entire compactified complexified Minkowski space
$\m$, only a codimension two submanifold is left out, since 
$\m \setminus \left( \m^{\rm I} \cup \m^{\rm J} \right) = C(\ell) \cap C(\ell')$.

The intersection $\m^{\rm IJ}=\m^{\rm I} \cap \m^{\rm J}$ is given by
the set of all $x\in\m^{\rm I}$ such that $\det(X)=x_{11'}x_{22'}-x_{12'}x_{21'}\neq0$.
Equivalently, this is also the set of all $y\in\m^{\rm J}$ such that 
$\det(Y)=y_{11'}y_{22'}-y_{12'}y_{21'}\neq0$. The gluing map 
$\tau:\m^{\rm IJ}\seta\m^{\rm IJ}$ relates the local coordinates on
$\m^{\rm I}$ and $\m^{\rm J}$ in the following way:
\begin{eqnarray*}
x_{11'} = \frac{y_{22'}}{\det(Y)} & \ \ \ &
x_{12'} = -\frac{y_{12'}}{\det(Y)} \\
x_{21'} = -\frac{y_{21'}}{\det(Y)} & \ \ \ &
x_{22'} = \frac{y_{11'}}{\det(Y)}
\end{eqnarray*}

Since $\m$ is an algebraic variety, it can also be 
characterized via its homogeneous coordinate algebra:
$$ \cm = \cpx[D,D',z_{11'},z_{12'},z_{21'},z_{22'}]_h/{\cal I}_m $$
where ${\cal I}_m$ is the ideal generated by the quadric (\ref{quad}). 
The subscript ``$h$" means that $\cm$ consists only of 
the homogeneous polynomials. 

In this picture, the coordinate algebras of the affine varieties
$\m^{\rm I}$ and $\m^{\rm J}$ introduced above can be interpreted
as certain localizations of the quadratic algebra $\cm$. Indeed, the
coordinate rings for $\m^{\rm I}$ and $\m^{\rm J}$ are respectively
given by:
\begin{eqnarray*}
\cm^{\rm I} = \cm[D^{-1}]_0 = 
\cpx\left[ x_{11'},x_{12'},x_{21'},x_{22'} \right] \\
\cm^{\rm J} = \cm[D'^{-1}]_0 = 
\cpx\left[ y_{11'},y_{12'},y_{21'},y_{22'} \right]
\end{eqnarray*}  
where the subscript ``$0$" means that we take only the degree zero part of
the localized graded algebra. 

Finally, we note that $\cm^{\rm I}$ and $\cm^{\rm J}$ can be made
isomorphic by adjoining the inverses of the determinants $\det(X)$
and $\det(Y)$, respectively. Indeed, define the matrices of generators:
\begin{equation} \label{XY} \begin{array}{ccc}
X = \left( \begin{array}{cc}
x_{11'} & x_{12'} \\ x_{21'} & x_{22'}
\end{array} \right) & \ \ &
Y = \left( \begin{array}{cc}
y_{22'} & -y_{12'} \\ -y_{21'} & y_{11'}
\end{array} \right)
\end{array} \end{equation}
The map $\eta:\cm^{\rm I}[\det(X)^{-1}] \rightarrow \cm^{\rm J}[\det(Y)^{-1}]$
given by:
$$ \eta(X) = \frac{Y}{\det(Y)} $$
is an isomorphism. It is the algebraic analogue of the gluing map $\tau$
described above in geometric context, while
$\cm^{\rm I}[\det(X)^{-1}] \simeq \cm^{\rm J}[\det(Y)^{-1}]$ plays the
role of the intersection $\m^{\rm IJ}=\m^{\rm I}\cap\m^{\rm J}$.

%------------------------------------------------------------

\subsection{The quantum group $SU(2)_q$ and its extension}

Let $q$ be a formal parameter.
Recall that the quantum group $GL(2)_q$ is the bialgebra over $\cpx$
generated by $g_{11'},g_{12'},g_{21'},g_{22'}$ subject to the
following commutation relations, see e.g. \cite{CP}:
\begin{eqnarray}
g_{11'}g_{12'}=q^{-1}g_{12'}g_{11'} & \ \ \ 
& g_{11'}g_{21'}=q^{-1}g_{21'}g_{11'} \nonumber\\
g_{12'}g_{22'}=q^{-1} g_{22'}g_{12'} & \ \ \ 
& g_{21'}g_{22'}=q^{-1}g_{22'}g_{21'} \label{gl2q1}
\end{eqnarray}
$$ g_{12'}g_{21'}=g_{21'}g_{12'} $$
\begin{equation} \label{gl2q2}
g_{11'}g_{22'} - q^{-1}g_{12'}g_{21'} = g_{22'}g_{11'} - qg_{21'}g_{12'}
\end{equation}
The comultiplication $\Delta$ and the counit $\varepsilon$ are given by:
\begin{equation} \label{comult}
\Delta(g_{rs'}) = \sum_{k=1,2} g_{rk'} \otimes g_{ks'}, \ \ \ r,s=1,2
\end{equation}
\begin{equation} \label{conunit}
\varepsilon(g_{rs'}) = \delta_{rs'}, \ \ \ r,s=1,2
\end{equation}

The expression (\ref{gl2q2}) is called the {\em quantum determinant} and
it is denoted by $\det_q(g)$. One easily checks that:
\begin{equation} \label{qdet.commutes}
g_{rs'}{\rm det}_q(g) = {\rm det}_q(g)g_{rs'}, \ \ \ r,s=1,2
\end{equation}
\begin{equation} \label{qdet.comult}
\Delta({\rm det}_q(g)) = {\rm det}_q(g)\otimes {\rm det}_q(g)
\end{equation}

In order to define a Hopf algebra structure on $GL(2)_q$ we adjoin the
inverse of the quantum determinant $\det_q(g)^{-1}$ to the generators
$g_{rs'}$ with the obvious relations. Then the antipode exists, and it
is given by:
$$ \begin{array} {ccc}
\gamma(g_{11'}) = g_{22'}{\rm det}_q(g)^{-1} & &
\gamma(g_{12'}) = -qg_{12'}{\rm det}_q(g)^{-1} \\
\gamma(g_{21'}) = -q^{-1}g_{21'}{\rm det}_q(g)^{-1} & &
\gamma(g_{22'}) = g_{11'}{\rm det}_q(g)^{-1} 
\end{array} $$
We will assume that the quantum group $GL(2)_q$ contains $\det_q(g)^{-1}$
and is therefore a Hopf algebra. 

The quantum group $SL(2)_q$ is defined by the quotient:
$$ SL(2)_q = GL(2)_q/\langle {\rm det}_q(g)=1 \rangle $$
It is a well defined Hopf algebra by 
equations (\ref{qdet.commutes}) and (\ref{qdet.comult}).

It is also useful to re-express the commutation relations for the quantum
group $GL(2)_q$ by means of the $R$-matrix; let
\begin{equation} \label{rmatrix}
R_{12} = \left( \begin{array}{cccc}
p^{-1} & 0 & 0 & 0 \\ 0 & 1 & p^{-1}-q & 0 \\
0 & p^{-1}-q^{-1} & 1 & 0 \\ 0 & 0 & 0 & p^{-1}
\end{array} \right), \ \ p=q^{\pm 1} \end{equation}
and let $T$ be the matrix of generators, i.e.:
$$ T = \left( \begin{array}{cc}
g_{11'} & g_{12'} \\ g_{21'} & g_{22'}
\end{array} \right) $$
Then the relations (\ref{gl2q1}) and (\ref{gl2q2}) can then 
be put in the following compact form:
$$ R_{12} T_1T_2 = T_2T_1 R_{12} $$
where $T_1=T\otimes\id$ and $T_2=\id\otimes T$ \cite{FRT}. Note that the 
commutation relations (\ref{gl2q1}) and (\ref{gl2q2})
do not depend on the parameter $p$. The $R$-matrix (\ref{rmatrix})
also satisfies the Hecke relation:
\begin{equation} \label{hecke1}
R_{12} - (R_{21})^{-1} = (p^{-1} - p) P
\end{equation}
where $R_{21}=R_{12}^{\rm t}$ (with the superscript ``t'' 
meaning transposition) and $P$ is the permutation matrix.
%$$ P = \left( \begin{array}{cccc}
%1 & 0 & 0 & 0 \\ 0 & 0 & 1 & 0 \\
%0 & 1 & 0 & 0 \\ 0 & 0 & 0 & 1
%\end{array} \right) $$
Equivalently, we also have for $\hat{R}_{12}=PR_{12}$: 
\begin{equation} \label{hecke2}
(\hat{R}_{12})^2 = (p^{-1} - p) \hat{R}_{12} + \id
\end{equation}

Next, we will extend the quantum groups $GL(2)_q$ and $SL(2)_q$  by
introducing a new generator $\delta$ and its inverse, satisfying the
following commutation relations with the quantum group generators:
\begin{equation} \label{deltacommut}
\begin{array} {ccc}
\delta g_{11'} = g_{11'} \delta & \ \ \  &
\delta g_{12'} = qg_{12'} \delta \\
\delta g_{21'} = q^{-1}g_{21'} \delta & \ \ \  &
\delta g_{22'} = g_{22'} \delta
\end{array} \end{equation}
In matrix form, the above relations become:
$$ \delta T Q^2 = Q^2 T \delta $$
where $Q$ is the following matrix:
\begin{equation} \label{matrixq}
Q = \left( \begin{array}{cc}
q^{\frac{1}{4}} & 0 \\ 0 & q^{-\frac{1}{4}}
\end{array} \right) \end{equation}
Strictly speaking, one should consider $q^{\frac{1}{4}}$ to be the
formal parameter, instead of $q$. However, fractional powers of $q$ will
rarely appear in this paper.

In other words, we define
$$ \widetilde{GL}(2)_q = GL(2)_q[\delta,\delta^{-1}]/(\ref{deltacommut})
\ \ \ {\rm and}\ \ \ 
\widetilde{SL}(2)_q = SL(2)_q[\delta,\delta^{-1}]/(\ref{deltacommut}) $$
The comultiplication, counit and antipode in $\widetilde{GL}(2)_q$ and
$\widetilde{SL}(2)_q$ are given by:
$$ \Delta(\delta) = \delta\otimes\delta, \ \ \
\varepsilon(\delta) = 1, \ \ \ 
\gamma(\delta)=\delta^{-1} $$
It is easy to check that $\widetilde{GL}(2)_q$ and
$\widetilde{SL}(2)_q$ 
satisfy the axioms of a Hopf algebra, and can thus be thought as quantum groups. 
Besides,
we have the identity:
$$ \gamma^2(g)=\delta^2 g \delta^{-2}, \ \ \
{\rm for\ all}\ g\in GL(2)_q, \ SL(2)_q $$
Thus conjugation by $\delta$ can be viewed as a square root of
the antipode squared.

In this paper, we will be primarily interested in the quantum group
$SL(2)_q$ and its extension $\widetilde{SL}(2)_q$. For these quantum
groups we define an involution $\dagger$ which fixes the formal parameter
(i.e. $q^\dagger = q$) and acts on the generators as follows:
\begin{equation} \label{involution}
\begin{array} {ccc}
g_{11'}^\dagger  = g_{22'}, & & g_{12'}^\dagger = -g_{21'} \\
g_{21'}^\dagger = -g_{12'}, & & g_{22'}^\dagger = g_{11'} \\
\delta^\dagger = \delta, & & (\delta^{-1})^\dagger = \delta^{-1}
\\
\end{array} \end{equation}
We extend the $\dagger$-involution to the quantum group $\widetilde{SL}(2)_q$
by requiring it to be a conjugate linear anti-homomorphism, that is:
\begin{equation} \label{ext.star} \begin{array}{l}
(xy)^\dagger = y^\dagger x^\dagger ,\ \ {\rm where}\ \  x,y\in\widetilde{SL}(2)_q;\\
a^\dagger=\overline{a},\ \ \ {\rm for\ all}\ \ a\in\cpx.
\end{array} \end{equation}

We define $SU(2)_q$ and $\widetilde{SU}(2)_q$ as the 
quantum groups $SL(2)_q$ and $\widetilde{SL}(2)_q$, respectively,
equipped with the $\dagger$-involution:
$$ SU(2)_q = (SL(2)_q,\dagger) 
\ \ \ {\rm and}\ \ \
\widetilde{SU}(2)_q=(\widetilde{SL}(2)_q,\dagger) $$

Note that instead of the formal parameter $q$ in the definitions
of the quantum groups, we could have used a positive real
number, which is automaticaly fixed by the involution $\dagger$.
However, we prefer to consider various specializations of formal $q$
later on, the most interesting one being the specialization to a
root of unity (see also Section \ref{ru} below).

%------------------------------------------------------------

\subsection{Quantum Minkowski space-time}

Let us introduce two new sets of variables on the extended
quantum group $\widetilde{SL}(2)_q$:
\begin{equation} \label{xvar} \begin{array} {lcr}
x_{11'} = \delta g_{11'} = g_{11'}\delta, & \ \ \ & 
x_{12'} = q^{-1/2}\delta g_{12'} = q^{1/2}g_{12'}\delta, \\
x_{21'} = q^{1/2}\delta g_{21'} = q^{-1/2}g_{21'}\delta, 
&\ \ \  & x_{22'} = \delta g_{22'} = g_{22'}\delta
\end{array} \end{equation}
and
\begin{equation} \label{yvar} \begin{array} {lcr}
y_{11'} = \delta^{-1} g_{11'} = g_{11'}\delta^{-1}, & \ \ \ & 
y_{12'} = q^{1/2}\delta^{-1} g_{12'} = q^{-1/2}g_{12'}\delta^{-1}, \\
y_{21'} = q^{-1/2}\delta^{-1} g_{21'} = q^{1/2}g_{21'}\delta^{-1}, 
& \ \ \  & y_{22'} = \delta^{-1} g_{22'} = g_{22'}\delta^{-1}
\end{array} \end{equation}
It is easy to determine their commutation relations:
\begin{equation} \label{xcommut1}
x_{11'}x_{12'} = x_{12'}x_{11'}, \ \ \ x_{21'}x_{22'} = x_{22'}x_{21'},
\end{equation}
\begin{equation} \label{xcommut3}
\left[ x_{11'} ,x_{22'} \right] + [x_{21'},x_{12'}] = 0
\end{equation}
\begin{equation} \label{xcommut2} \begin{array} {c}
x_{11'}x_{21'} = q^{-2}x_{21'}x_{11'}, \ \ \ x_{12'}x_{22'} = q^{-2}x_{22'}x_{12'}, \\
x_{21'}x_{12'} = q^2 x_{12'} x_{21'}
\end{array} \end{equation}
and
\begin{equation} \label{ycommut1}
y_{11'}y_{21'} = y_{21'}y_{11'}, \ \ \ y_{12'}y_{22'} = y_{22'}y_{12'}
\end{equation}
\begin{equation} \label{ycommut3}
\left[ y_{11'},y_{22'} \right] + [y_{12'},y_{21'}] = 0
\end{equation}
\begin{equation} \label{ycommut2} \begin{array} {c}
y_{11'}y_{12'} = q^{-2}y_{12'}y_{11'}, \ \ \ y_{21'}y_{22'} = q^{-2}y_{22'}y_{21'}, \\
y_{12'}y_{21'} = q^2 y_{21'} y_{12'}
\end{array} \end{equation}
The determinant condition $\det_q(g)=1$ in the new variables becomes:
\begin{equation} \label{xdtx} 
x_{11'}x_{22'} - x_{12'}x_{21'} = 
x_{22'}x_{11'} - x_{21'}x_{12'} = \delta^2
\end{equation}
and
\begin{equation} \label{ydty}
y_{11'}y_{22'} - y_{21'}y_{12'} =
 y_{22'}y_{11'} - y_{12'}y_{21'} = \delta^{-2}
 \end{equation}
Furthermore, the $\dagger$ involution acts on these new variables in the
following way:
\begin{equation} \label{star.xy}
\begin{array}{ccccccc}
x_{11'}^\dagger = x_{22'} & \ \ & x_{12'}^\dagger =-x_{21'} & \ \ &
x_{21'}^\dagger =-x_{12'} & \ \ & x_{22'}^\dagger = x_{11'} \\
y_{11'}^\dagger = y_{22'} & \ \ & y_{12'}^\dagger =-y_{21'} & \ \ &
y_{21'}^\dagger =-y_{12'} & \ \ & y_{22'}^\dagger = y_{11'}
\end{array}\end{equation}

We will consider the following subalgebras of $\widetilde{SL}(2)_q$:
$$
\cm^{\rm I}_q = \cpx[x_{11'},x_{12'},x_{21'},x_{22'}]/(\ref{xcommut1}-\ref{xcommut2})
$$
and
$$
\cm^{\rm J}_q = \cpx[y_{11'},y_{12'},y_{21'},y_{22'}]/(\ref{ycommut1}-\ref{ycommut2})
$$
regarding them as $q$-deformations of $\cm^{\rm I}$ and $\cm^{\rm J}$, the
coordinate algebras corresponding to the affine subspaces $\m^{\rm I}$ and 
$\m^{\rm J}$ of the compactified complexified Minkowski space-time $\m$.
Moreover, we introduce the $\dagger$-algebras:
$$ S^{\rm I}_q=(\cm^{\rm I}_q,\dagger)
\ \ \ {\rm and}\ \ \ 
S^{\rm J}_q=(\cm^{\rm J}_q,\dagger) $$
regarding them as $q$-deformations of
$S^4\setminus\{\infty\}$ and $S^4\setminus\{0\}$, respectively.

Note also that $\cm^{\rm I}_q$ and $\cm^{\rm J}_q$ are isomorphic (as Hopf algebras) 
once the generators $\delta$ and $\delta^{-1}$ are adjoined. Let
$X$ and $Y$ be again the matrices of generators as defined in (\ref{XY}); 
the map:
\begin{eqnarray}
\nonumber \eta: \cm^{\rm I}_q[\delta^{-1}] & \longrightarrow & \cm^{\rm J}_q[\delta] \\
\label{nc.glueing} X & \mapsto & \delta^2 Y 
\end{eqnarray}
is an isomorphism. More explicitly, in coordinates:
\begin{equation} \begin{array}{ccc}
x_{11'} \mapsto \delta^2 y_{22'} & \ \ \ & 
x_{12'} \mapsto -\delta^2 y_{12'} \\
x_{21'} \mapsto -\delta^2 y_{21'} & \ \ \ & 
x_{22'} \mapsto \delta^2 y_{11'}
\end{array} \end{equation}

The inverse map $\eta^{-1}$ is given by $Y\mapsto X \delta^{-2}$. 
Moreover, $\eta(X^\dagger) = \eta(X)^\dagger$, so that
$\eta(f^\dagger)=\eta(f)^\dagger$ for all $f\in\cm^{\rm I}_q[\delta^{-1}]$.

Geometrically, the algebras $\cm^{\rm I}_q[\delta^{-1}]$ and
$\cm^{\rm J}_q[\delta]$ play the role of the intersection 
$\m^{\rm IJ}=\m^{\rm I}\cap\m^{\rm J}$, while $\eta$ plays the
role of the gluing map $\tau$.

Furthermore, as in the case of the quantum group $GL(2)_q$, the 
above commutation relations can also be put in a compact form by
means of an $R$-matrix. It is now convenient to consider the
following matrix of generators:
\begin{equation} \label{XY2} \begin{array}{ccc}
X = \left( \begin{array}{cc}
x_{11'} & x_{12'} \\ x_{21'} & x_{22'}
\end{array} \right) 
& \ \ \ &
Y = \left( \begin{array}{cc}
y_{11'} & y_{12'} \\ y_{21'} & y_{22'}
\end{array} \right)
\end{array} \end{equation}
Notice that:
$$ \begin{array}{ccc}
X = Q^{-1} \delta T Q = Q T \delta Q^{-1}
& \ \ \ {\rm and}\ \ \ &
Y = Q \delta^{-1} T Q^{-1} = Q^{-1} T \delta^{-1} Q
\end{array} $$
where $Q$ was defined in (\ref{matrixq}).
In order to write down the above commutation relations in matrix form,
define $X_1=X\otimes\id$, $X_2=\id\otimes X$ and similarly $Q_1=Q\otimes\id$, 
$Q_2 = \id\otimes Q$ . Therefore:
\begin{eqnarray*}
X_1X_2 & = & 
\left( Q_1^{-1} \delta T_1 Q_1 \right) \left( Q_2 T_2 \delta Q_2^{-1} \right) = 
\left( Q_1^{-1} Q_2 \right) \delta T_1T_2 \delta \left( Q_1Q_2^{-1} \right) = \\
& = & \left( Q_1^{-1} Q_2 R^{-1}_{12} Q_2 Q_1^{-1} \right) X_2X_1 
\left( Q_2^{-1} Q_1 R_{12} Q_1 Q_2^{-1} \right)
\end{eqnarray*}
Thus we define $R^{\rm I}_{12} = Q_2^{-1} Q_1 R_{12} Q_1 Q_2^{-1}$; 
more precisely:
\begin{equation} \label{rmatrixI}
R^{\rm I}_{12} = \left( \begin{array}{cccc}
p^{-1} & 0 & 0 & 0 \\ 0 & q^{-1} & p^{-1}-q & 0 \\
0 & p^{-1}-q^{-1} & q & 0 \\ 0 & 0 & 0 & p^{-1}
\end{array} \right), \ \ p=q^{\pm1}
\end{equation}
The commutation relations (\ref{xcommut1}-\ref{xcommut2}) can then be
presented in matrix form:
$$ R^{\rm I}_{12} X_1X_2 =  X_2X_1 R^{\rm I}_{12} $$

Performing a similar calculation for the $y_{rs'}$ variables, 
we obtain:
\begin{equation} \label{rmatrixJ}
R^{\rm J}_{12} = Q_2 Q_1^{-1} R_{12} Q_1^{-1} Q_2 = \left( \begin{array}{cccc}
p^{-1} & 0 & 0 & 0 \\ 0 & q & p^{-1}-q & 0 \\
0 & p^{-1}-q^{-1} & q^{-1} & 0 \\ 0 & 0 & 0 & p^{-1}
\end{array} \right),  \ \ p=q^{\pm1}
\end{equation}
with the commutation relations (\ref{ycommut1}-\ref{ycommut2})
being given by:
$$ R^{\rm J}_{12} Y_1Y_2 =  Y_2Y_1 R^{\rm J}_{12} $$

%---------------------------------------------------------------
%---------------------------------------------------------------

\section{Differential forms on quantum Minkowski space-time} \label{forms}

\subsection{Differential forms on Minkowski space-time}

Since $\m^{\rm I}$ and $\m^{\rm J}$ are affine spaces, 
their modules of differential forms are very simple to describe. Indeed, 
recall that:
$$ \cm^{\rm I}=\cpx[x_{11'},x_{12'},x_{21'},x_{22'}] $$
Therefore, the module of differential 1-forms is given by the free
$\cm^{\rm I}$-module generated by $dx_{rs'}$:
$$ \Omega^1_{\cm^{\rm I}} = \cm^{\rm I} \langle dx_{rs'} \rangle $$
while the module of differential 2-forms is given by the free
$\cm^{\rm I}$-module generated by $dx_{rs'}\wedge dx_{kl'}$:
$$ \Omega^2_{\cm^{\rm I}} = \Lambda^2(\Omega^1_{\cm^{\rm I}}) = 
\cm^{\rm I} \langle dx_{rs'}\wedge dx_{kl'} \rangle $$
with $r,s,k,l=1,2$.

The action of the de Rham operator $d: \cm^{\rm I}\seta\Omega^1_{\cm^{\rm I}}$
is given on the generators as $x_{rs'}\mapsto dx_{rs'}$,
and it is then extended to the whole $\cm^{\rm I}$ by $\cpx$-linearity and 
the Leibnitz rule:
\begin{equation} \label{l}
d(fg) = g df + f dg
\end{equation}
where $f,g \in \cm^{\rm I}$. One also defines the de Rham operator
$d : \Omega^1_{\cm^{\rm I}}\seta\Omega^2_{\cm^{\rm I}}$ on the
generators as $f dx_{rs'}\mapsto df \wedge dx_{rs'}$, also
extending it by $\cpx$-linearity and the Leibnitz rule (\ref{l}).

The modules of differential forms and de Rham operators over $\cm^{\rm J}$
are similarly described. 

Now let $\Omega^2_{E^4}$ denote the bundle of 2-forms on Euclidean
space $E^4$ with coordinates $\{x^1,x^2,x^3,x^4\}$.
Recall that the Hodge operator $*:\Omega^2_{E^4}\seta\Omega^2_{E^4}$
is defined as follows:
$$ \begin{array}{ccc}
*dx^1\wedge dx^2 = dx^3\wedge dx^4 & \ \ \ \ &
*dx^3\wedge dx^4 = dx^1\wedge dx^2 
\end{array} $$
$$ \begin{array}{ccc}
*dx^1\wedge dx^3 = -dx^2\wedge dx^4 & \ \ \ \ &
*dx^2\wedge dx^4 = -dx^1\wedge dx^3 
\end{array} $$
$$ \begin{array}{ccc}
*dx^1\wedge dx^4 = dx^2\wedge dx^3 & \ \ \ \ &
*dx^2\wedge dx^3 = dx^1\wedge dx^4 
\end{array} $$
We can then use the relation between Euclidean and twistor coordinates 
on $\m^{\rm I}=E^4\otimes\cpx$ given by (\ref{coords1}) to
express the action of the Hodge operator on $\Omega^2_{\cm^{\rm I}}$. One obtains:
$$ \begin{array}{ccc}
*dx_{11'}\wedge dx_{12'} = dx_{11'}\wedge dx_{12'} & \ \ \ \ &
*dx_{11'}\wedge dx_{21'} = - dx_{11'}\wedge dx_{21'}
\end{array} $$
$$ \begin{array}{ccc}
*dx_{11'}\wedge dx_{22'} = - dx_{12'}\wedge dx_{21'} & \ \ \ \ &
*dx_{12'}\wedge dx_{21'} = - dx_{11'}\wedge dx_{22'}
\end{array} $$
$$ \begin{array}{ccc}
*dx_{12'}\wedge dx_{22'} = - dx_{12'}\wedge dx_{22'} & \ \ \ \ &
*dx_{21'}\wedge dx_{22'} = dx_{21'}\wedge dx_{22'}
\end{array} $$

Clearly $*^2=1$, thus the complexified Hodge operator 
$*:\Omega^2_{\cm^{\rm I}}\seta\Omega^2_{\cm^{\rm I}}$
induces a splitting of $\Omega^2_{\cm^{\rm I}}$ into two
submodules corresponding to eigenvalues $\pm1$. More explicitly, we have:
$$ \Omega^{2}_{\cm^{\rm I}}=\Omega^{2,+}_{\cm^{\rm I}}\oplus
\Omega^{2,-}_{\cm^{\rm I}} $$
where
$$ \Omega^{2,+}_{\cm^{\rm I}} = \cm^{\rm I} \langle dx_{11'}\wedge dx_{12'}, 
dx_{21'}\wedge dx_{22'}, dx_{11'}\wedge dx_{22'} - dx_{12'}\wedge dx_{21'}
\rangle $$
$$
\Omega^{2,-}_{\cm^{\rm I}} = \cm^{\rm I} \langle dx_{11'}\wedge dx_{21'}, 
dx_{12'}\wedge dx_{22'}, dx_{11'}\wedge dx_{22'} + dx_{12'}\wedge dx_{21'}
\rangle $$

\paragraph{Connection and curvature}
Let $E$ be a $\cm^{\rm I}$-module;
a {\em connection} on $E$ is a $\cpx$-linear map:
$$ \nabla: E \seta E\otimes_{\cm^{\rm I}}\Omega^1_{\cm^{\rm I}} $$
satisfying the Leibnitz rule:
\begin{equation} \label{leib1}
\nabla(f\sigma)=\sigma\otimes df + f\nabla\sigma
\end{equation}
where $f\in\cm^{\rm I}$ and $\sigma\in E$.
The connection $\nabla$ also acts on 1-differentials, being defined as the
additive map:
$$ \nabla: E\otimes_{\cm^{\rm I}}\Omega^1_{\cm^{\rm I}} \seta 
E\otimes_{\cm^{\rm I}}\Omega^2_{\cm^{\rm I}} $$
satisfying:
\begin{equation} \label{leib2}
\nabla(\sigma\otimes\omega)=\sigma\otimes d\omega +
\omega\wedge\nabla\sigma
\end{equation}
where $\omega\in\Omega^1_{\cm^{\rm I}}$.

Moreover, two connections $\nabla$ and $\nabla'$ are said to be gauge equivalent
if there is $g\in{\rm Aut}_{\cm^{\rm I}}(E)$ such that
$\nabla = g^{-1} \nabla' g$.

The curvature $F_\nabla$ is defined by the composition:
$$ E \stackrel{\nabla}{\longrightarrow} 
E\otimes_{\cm^{\rm I}}\Omega^1_{\cm^{\rm I}}
\stackrel{\nabla}{\longrightarrow} E\otimes_{\cm^{\rm I}}
\Omega^2_{\cm^{\rm I}} $$
and it is easy to check that it is actually $\cm^{\rm I}$-linear. 
Therefore, $F_\nabla$ can be regarded as an element of
${\rm End}_{\cm^{\rm I}}(E)\otimes_{\cm^{\rm I}}\Omega^2_{\cm^{\rm I}}$.
Furthermore, if $\nabla$ and $\nabla'$ are gauge equivalent, then
there is $g\in{\rm Aut}_{\cm^{\rm I}}(E)$ such that
$F_\nabla = g^{-1} F_{\nabla'} g$.

If $E$ is projective (hence free), any connection $\nabla$ can be encoded
into a matrix 
$A\in{\rm End}_{\cm^{\rm I}}(E)\otimes_{\cm^{\rm I}}\Omega^1_{\cm^{\rm I}}$. 
Simply choose a basis  $\{\sigma_k\}$ for $E$, and let 
$s=\sum a^k\sigma_k$,
so that:
$$ \nabla s = \sum_k \left( \sigma_k\otimes da^k + a^k\nabla\sigma_k \right) $$
Thus it is enough to know how $\nabla$ acts on the basis $\{\sigma_k\}$:
$$ \nabla\sigma_k = \sum_{l,\alpha} A_k^{l,\alpha}\sigma_l\otimes dx_\alpha $$
and we define $A$ as the matrix with entries given by the 1-forms
$\sum_\alpha A_k^{l,\alpha}dx_\alpha$.

Conversely, given 
$A\in{\rm End}_{\cm^{\rm I}}(E)\otimes_{\cm^{\rm I}}\Omega^1_{\cm^{\rm I}}$,
we define the connection:
$$ \nabla_A s = \sum_{k,l,\alpha}
\left( \sigma_k\otimes da^k + A_k^{l,\alpha}a^k\sigma_l\otimes dx_\alpha \right) $$

%------------------------

\subsection{Differential forms on the quantum group $SU(2)_q$}

Let us now recall a few facts regarding the exterior algebra over 
the relevant quantum groups \cite{S,WZ,W}.
The module of 1-forms over the quantum group $GL(2)_q$, which we 
shall denote by $\Omega^1_{GL}$, is the $GL(2)_q$-bimodule 
generated by $dg_{rs'}$ satisfying the following relations 
(written in matrix form):
\begin{equation} \label{nc1forms}
R_{12} T_1  dT_2 = dT_2  T_1 (R_{21})^{-1}
\end{equation}
where $R_{12}$ is again the matrix (\ref{rmatrix}) 
and $dT_2 = \id\otimes dT$, with:
$$ dT = \left( \begin{array}{cc}
dg_{11'} & dg_{12'} \\ dg_{21'} & dg_{22'} 
\end{array} \right) $$
Similarly, the module of 2-forms $\Omega^2_{GL}$, is the $GL(2)_q$-bimodule 
generated by $dg_{rs'}\wedge dg_{kl'}$, which satisfy the relations 
(written in matrix form):
\begin{equation} \label{nc2forms}
R_{12} dT_1 \wedge dT_2 = - dT_2 \wedge dT_1 (R_{21})^{-1}
\end{equation}
where $dT_1 = dT \otimes\id$.
Furthermore, the commutation relations between $g_{mn'}$ and  
$dg_{rs'}\wedge dg_{kl'}$ can be deduced from (\ref{nc1forms})
and (\ref{nc2forms}) as follows. Let $dT_3 = \id\otimes\id\otimes dT$.
Denoting $R_{ba}=(R_{ab})^{\rm t}$, we have 
$R_{12} T_1  dT_2 R_{21} = dT_2  T_1$ and
$R_{13} T_1  dT_3 R_{31} = dT_3  T_1$. Therefore,
\begin{eqnarray*}
dT_3 \wedge dT_2 T_1 & = & 
dT_3 \wedge (R_{12} T_1 dT_2 R_{21}) = R_{12} dT_3 T_1 \wedge
dT_2 R_{21} = \\
& = & R_{12}R_{13} T_1  dT_3 R_{31} \wedge dT_2 R_{21} = \left( R_{12}R_{13} \right) T_1  
dT_3 \wedge dT_2 \left( R_{12}R_{13} \right)^{\rm t}
\end{eqnarray*} 

The noncommutative de Rham operators are given by their
action on the generators as follows:  
$$ \begin{array}{ccc}
d : GL(2)_q \seta \Omega^1_{GL} & \ \ \ \ & 
d : \Omega^1_{GL} \seta \Omega^2_{GL} \\
g_{rs'} \mapsto dg_{rs'} & \ \ \ \ & 
g_{kl'}dg_{rs'} \mapsto dg_{kl'} \wedge dg_{rs'} 
\end{array} $$
This is then extended to the whole $GL(2)_q$ and $\Omega^1_{GL}$
by $\cpx$-linearity and the Leibnitz rule:
$$ d(f_1f_2) =  df_1 f_2  + f_1 df_2 $$
for all $f_1,f_2\in GL(2)_q$. 

The modules $\Omega^1_{GL}$ and $\Omega^2_{GL}$ also have natural
involutions, extended from $\dagger$ in a natural way, namely: 
$$ (f dg_{rs'})^\dagger = dg_{rs'}^\dagger f^\dagger $$
$$ (f dg_{rs'} \wedge dg_{kl'})^\dagger = -
dg_{kl'}^\dagger\wedge dg_{rs'}^\dagger f^\dagger $$

To get the modules of forms on $SL(2)_q$ it is enough to take the 
quotient by the appropriate relations:
$$ \Omega^1_{SL} = \Omega^1_{GL} \left/ 
\begin{array}{c} \det_q T = 1 \\ d (\det_q T) = 0 \end{array} 
\right. \ \ \ {\rm and} \ \ \ \Omega^2_{SL}=\Lambda^2(\Omega^1_{SL}) $$
Finally, modules of forms on $SU(2)_q$ are then defined as the pairs:
$$ \Omega^k_{SU} = ( \Omega^k_{SL} , \dagger ), \ \ \ k=1,2 $$

It is also important for our purposes to describe the modules of 
1- and 2-forms on the extended quantum group $\widetilde{SL}(2)_q$. 
To do that, we add the generators $\delta$ and  $d\delta$ satisfying
the relations (written in matrix form): 
$$ \begin{array}{ccc}
\delta dT Q^2 = Q^2 dT \delta
&  \ \ {\rm and} \ \ &
d\delta dT Q^2 = - Q^2 dT d\delta
\end{array} $$

We extend the involution $\dagger$ to $\Omega^k_{\widetilde{SL}}$ by
declaring that $(d\delta)^\dagger = d\delta$; thus we have:
$$ \Omega^k_{\widetilde{SU}} = ( \Omega^k_{\widetilde{SL}} , \dagger ), \ \ \ k=1,2 $$

Let us now introduce the noncommutative analogue 
of the Hodge operator on $\Omega^2_{GL}$. This will serve as a model 
on our Definition of self-dual and anti-self-dual 2-forms on 
quantum space-time.

Recall that $R_{12}$ denotes the $R$-matrix (\ref{rmatrix}) for $GL(2)_q$. Define:
\begin{equation} \label{p+p-}
P^+ = \frac{\hat{R}_{12} + p \id}{p+p^{-1}}
\ \ \ {\rm and} \ \ \
P^- = \frac{-\hat{R}_{12} + p^{-1} \id}{p+p^{-1}}
\end{equation}
Clearly, $P^+ + P^- = \id$. Moreover, using the Hecke relations (\ref{hecke2}), 
one easily checks that $(P^+)^2 = P^+$ and $(P^-)^2 = P^-$. 
Now the Hodge $*$-operator is defined on the generators of $\Omega^2_{GL}$ by:
\begin{equation} \label{hodgeGL}
* dT_1 \wedge dT_2 = \left( P^+ - P^- \right) dT_1 \wedge dT_2 
\end{equation}

%-----------------------------------------------------

\subsection{Differential forms on quantum Minkowski space-time}

In analogy with the classical case, we define the module of 1-forms over
the algebra $\cm^{\rm I}_q$, which we shall denote by $\Omega^1_{\cm^{\rm I}_q}$,
as the $\cm^{\rm I}_q$-bimodule generated by:
$$ dX =  
\left( \begin{array}{cc}
dx_{11'} & dx_{12'} \\ dx_{21'} & dx_{22'} 
\end{array} \right) = Q^{-1} \delta dT Q = Q dT \delta Q^{-1} $$

Moreover, the generators $dx_{rs'}$ satisfy the following relations 
(written in matrix form):
\begin{eqnarray*} 
X_1  dX_2 & = & 
\left( Q_1^{-1}Q_2 \right) \delta T_1 dT_2 \delta \left( Q_1Q_2^{-1} \right) = 
\left( Q_1^{-1}Q_2 \right) \delta R_{12}^{-1} dT_2 T_1 R_{21}^{-1} \delta \left( Q_1 Q_2^{-1} \right) = \\
& = & \left( Q_1^{-1}Q_2 R_{12}^{-1} Q_2 Q_1^{-1} \right) dX_2 X_1 \left( Q_2^{-1} Q_1 R_{21}^{-1} Q_1 Q_2^{-1} \right)  
\end{eqnarray*} 
where $R_{12}$ is again the matrix (\ref{rmatrix}) and $dX_2 = \id\otimes dX$.
Thus using the $R$-matrix for $\cm^{\rm I}_q$ (\ref{rmatrixI}) and defining 
$R^{\rm I}_{21}=Q_1^{-1}Q_2R_{21}Q_1^{-1}Q^{2}$, we obtain:
\begin{equation} \label{xnc1forms}
R^{\rm I}_{12} X_1 dX_2 = dX_2  X_1 (R^{\rm I}_{21})^{-1} 
\end{equation}

Similarly, the module of 2-forms $\Omega^2_{\cm^{\rm I}_q}$, is the 
$\cm^{\rm I}_q$-bimodule generated by $dx_{rs'}\wedge dx_{kl'}$ satisfying 
the relations below (written in matrix form), which can be deduced from 
(\ref{xnc1forms}):
\begin{equation} \label{xnc2forms}
R^{\rm I}_{12} dX_1 \wedge dX_2 = - dX_2 \wedge dX_1 (R^{\rm I}_{21})^{-1} 
\end{equation}
where $dX_1 = dX \otimes \id$.

Performing the same calculations for $\cm^{\rm J}_q$, we conclude that 
$\Omega^1_{\cm^{\rm J}_q}$ and $\Omega^2_{\cm^{\rm J}_q}$ are the 
$\cm^{\rm J}_q$-bimodules generated by $dy_{rs'}$ and $dy_{rs'}\wedge dy_{kl'}$, 
respectively, satisfying the following relations (written in matrix form):
\begin{equation} \label{yncforms} \begin{array}{ccc}
R^{\rm J}_{12} Y_1 dY_2 = dY_2  Y_1 (R^{\rm J}_{21})^{-1}
& \ \ {\rm and}\ \ &
R^{\rm J}_{12} dY_1 \wedge dY_2 = - dY_2 \wedge dY_1 (R^{\rm J}_{21})^{-1}
\end{array} \end{equation}
where $R^{\rm J}_{12}$ is the $R$-matrix for $\cm^{\rm J}_q$ (\ref{rmatrixJ}), 
$R^{\rm J}_{21}=Q_2^{-1}Q_1R_{21}q_1Q_2^{-1}$ and
$$ dY =  \left( \begin{array}{cc}
dy_{11'} & dy_{12'} \\ dy_{21'} & dy_{22'} 
\end{array} \right) = Q\delta^{-1} dT Q^{-1} = Q^{-1} dT \delta^{-1} Q $$
with $dY_1 = dY \otimes\id$ and $dY_2 = \id\otimes dY$. 

We now introduce the concept of anti-self-duality of quantum 2-forms
over the quantum Minkowski space-time. Notice that the matrices 
$R^{\rm I}_{12}$, $R^{\rm I}_{21}$ and $R^{\rm J}_{12}$, $R^{\rm J}_{21}$ 
also satisfy Hecke relations (recall that $P$ denotes the permutation matrix):
$$ R^{\rm I}_{12} - (R^{\rm I}_{21})^{-1} = (p^{-1}-p) P
\ \ \ {\rm and}\ \ \ 
R^{\rm J}_{12} - (R^{\rm J}_{21})^{-1} = (p^{-1}-p) P $$
Moreover, defining $\hat{R}^{\rm I}_{12}=PR^{\rm I}_{12}$
and $\hat{R}^{\rm J}_{12}=PR^{\rm J}_{12}$, we obtain:
\begin{equation} \label{heckeIJ}
(\hat{R}^{\rm I}_{12})^2 = \id + (p^{-1}-p) \hat{R}^{\rm I}_{12}
\ \ \ {\rm and}\ \ \ 
(\hat{R}^{\rm J}_{12})^2 = \id + (p^{-1}-p) \hat{R}^{\rm J}_{12}
\end{equation}
in analogy with (\ref{hecke2}). Therefore we can proceed as in the case of 
$\widetilde{GL}(2)_q$ discussed above and define the projectors:
$$ P^{{\rm I}+} = \frac{\hat{R}^{\rm I}_{12} + p \id}{p+p^{-1}}
\ \ \ {\rm and} \ \ \
P^{{\rm I}-}  = \frac{-\hat{R}^{\rm I}_{12} + p^{-1} \id}{p+p^{-1}} $$
$$ P^{{\rm J}+} = \frac{\hat{R}^{\rm J}_{12} + p \id}{p+p^{-1}}
\ \ \ {\rm and} \ \ \
P^{{\rm J}-}  = \frac{-\hat{R}^{\rm J}_{12} + p^{-1} \id}{p+p^{-1}} $$
Now we define the Hodge operator on $\Omega^2_{\cm^{\rm I}_q}$ 
(in matrix form):
$$ * dX_1 \wedge dX_2 =  \left( P^{{\rm I}+} - P^{{\rm I}-} \right) dX_1 \wedge dX_2 $$
Since $*^2=\id$, the module $\Omega^2_{\cm_q^{\rm I}}$ can be decomposed
into two submodules corresponding to eigenvalues $\pm1$. Denote such
submodules by $\Omega^{2,+}_{\cm_q^{\rm I}}$ and
$\Omega^{2,-}_{\cm_q^{\rm I}}$. 

In order to compare with the commutative case, it is instructive 
to write down their bases, which are given by the entries of the
matrices $P^{{\rm I}+} dX_1 \wedge dX_2$ and $P^{{\rm I}-} dX_1 \wedge dX_2$. 
After applying the commutation relations (\ref{nc2forms}), we conclude
that:
\begin{equation} \label{nc-sd2f} 
\Omega^{2,+}_{\cm^{\rm I}_q} = \cm^{\rm I}_q \langle dx_{11'}\wedge dx_{12'}, 
dx_{21'}\wedge dx_{22'}, dx_{11'}\wedge dx_{22'} - dx_{12'}\wedge dx_{21'}
\rangle
\end{equation}
\begin{equation} \label{nc-asd2f}
\Omega^{2,-}_{\cm^{\rm I}_q} = \cm^{\rm I}_q \langle dx_{11'}\wedge dx_{21'}, 
dx_{12'}\wedge dx_{22'}, dx_{11'}\wedge dx_{22'} + dx_{12'}\wedge dx_{21'}
\rangle
\end{equation}
in complete analogy with the commutative case.

Finally, replacing all the I's by J's, we define the Hodge operator on  
$\Omega^2_{\cm_q^{\rm J}}$ as well: 
$$ * dY_1 \wedge dY_2 = 
\left( P^{{\rm J}+} - P^{{\rm J}-} \right) dY_1 \wedge dY_2 $$

\paragraph{Connection and curvature}
Let $E$ be a right $\cm^{\rm I}_q$-module. In analogy with the commutative case,
a connection on $E$ is a $\cpx$-linear map:
$$ \nabla: E \seta E\otimes_{\cm^{\rm I}_q}\Omega^1_{\cm^{\rm I}_q} $$
satisfying the Leibnitz rule:
$$ \nabla(\sigma f)=\sigma\otimes df + \nabla(\sigma) f $$
where $f\in\cm^{\rm I}_q$ and $\sigma\in E$.
The connection $\nabla$ also acts on 1-forms, being defined as the
$\cpx$-linear map:
$$ \nabla: E\otimes_{\cm^{\rm I}_q}\Omega^1_{\cm^{\rm I}_q} \seta 
E\otimes_{\cm^{\rm I}_q}\Omega^2_{\cm^{\rm I}_q} $$
satisfying:
$$ \nabla(\sigma\otimes\omega)=\sigma\otimes d\omega +
\nabla\sigma \wedge \omega $$
where $\omega\in\Omega^1_{\cm^{\rm I}_q}$.

Moreover, two connections $\nabla$ and $\nabla'$ are said to be gauge equivalent
if there is $g\in{\rm Aut}_{\cm^{\rm I}_q}(E)$ such that
$\nabla = g^{-1} \nabla' g$.

The curvature $F_\nabla$ is defined by the composition:
$$ E \stackrel{\nabla}{\longrightarrow} 
E\otimes_{\cm^{\rm I}_q}\Omega^1_{\cm^{\rm I}_q}
\stackrel{\nabla}{\longrightarrow} E\otimes_{\cm^{\rm I}_q}
\Omega^2_{\cm^{\rm I}_q} $$
and it is easy to check that $F_\nabla$ is actually right $\cm^{\rm I}_q$-linear. 
Therefore, $F_\nabla$ can be regarded as an element of
${\rm End}_{\cm^{\rm I}_q}(E)\otimes_{\cm^{\rm I}_q}\Omega^2_{\cm^{\rm I}_q}$.
Furthermore, if $\nabla$ and $\nabla'$ are gauge equivalent, then
there is $g\in{\rm Aut}_{\cm^{\rm I}}(E)$ such that
$F_\nabla = g^{-1} F_{\nabla'} g$.
A connection $\nabla$ is said to be anti-self-dual if 
$F_\nabla\in
{\rm End}_{\cm^{\rm I}_q}(E)\otimes_{\cm^{\rm I}_q}\Omega^{2,-}_{\cm_q^{\rm I}}$.

Finally, if $E$ is projective (though not necessarily free in this case), 
any connection $\nabla$ can be encoded into the {\em connection matrix} 
$A\in{\rm End}_{\cm^{\rm I}_q}( \oplus^n\cm^{\rm I}_q )\otimes_{\cm^{\rm I}_q}\Omega^1_{\cm^{\rm I}_q}$
(where $n$ is the rank of $E$) in the following way. First, we recall the following basic result from algebra
(see e.g. \cite{Lm}):

\begin{theorem}[Dual basis theorem]
A finitely generated right $R$-module $M$ is projective of rank $n$ if and only if there
are elements $\sigma_k\in M$ and $\rho^k\in M^\vee={\rm Hom}_R(M,R)$ for $k=1,...,n$
such that $m=\sum_k \sigma_k \rho^k(m)$, for any $m\in M$.
\end{theorem}

Let $\{\sigma_k,\rho^k\}$ be a dual basis for $E$, so that any $s\in E$ 
can be written as $s=\sum_k \sigma_k \rho^k(s)$; applying the Leibnitz 
rule, we get:
$$ \nabla s = 
\sum_k \left( \sigma_k\otimes d(\rho^k(s)) + \nabla(\sigma_k) \rho^k(s) \right) $$
Thus, as in the commutative case, it is enough to know how $\nabla$ acts 
on $\{\sigma_k\}$; we set:
$$ \nabla\sigma_k = \sum_{l} \sigma_l A_k^{l},
\ \ \ {\rm with}\ \ \ A_k^l = \rho^l\otimes\id(\nabla\sigma_k) \in \Omega^1_{\cm^{\rm I}_q} $$
and we define $A$ as the matrix with entries given by the 1-forms
$A_k^l$. Clearly, $A$ depends on the choice of dual basis; changing the dual
basis amounts to a change of gauge for $A$.

Conversely given a matrix
$A\in{\rm End}_{\cm^{\rm I}_q}(\oplus^n\cm^{\rm I}_q)\otimes_{\cm^{\rm I}_q}\Omega^1_{\cm^{\rm I}_q}$,
we define the connection:
$$ \nabla_A s = \sum_{l,k} \sigma_k\otimes d(\rho^k(s)) + 
\sigma_l A_k^l \rho^k(s) $$

%---------------------------------------------------------------
%---------------------------------------------------------------

\section{Construction of quantum instantons} \label{ward}

\subsection{Classical instantons and the ADHM data}

As we mentioned in Introduction, solutions of the classical ASDYM
equations can be constructed from certain linear data, so-called {\em ADHM
data}. We will now briefly review some relevant facts regarding this
correspondence.

Recall that the ASDYM equation can be defined for a complex vector bundle
$E$ over any four dimensional Riemannian manifold $X$, provided with a 
connection $\nabla$. This connection is said to be {\em anti-self-dual}
if the corresponding curvature 2-form $F_\nabla$ satisfies the equation
($*$ is the Hodge operator on 2-forms):
$$ *F_\nabla = - F_\nabla $$
i.e., if $F_\nabla\in{\rm End}(E)\otimes\Omega^{2,-}_X$. An ASD connection
is usually called an {\em instanton} if the integral:
$$ c = \frac{1}{8\pi^2} \int_X {\rm Tr}(F_\nabla \wedge F_\nabla) $$
converges. If $X$ is a compact manifold, $c$ is actually an integer,
so-called {\em instanton number} or {\em charge}, and coincides with 
the second Chern class of $E$. In view of the symmetry between 
self-dual and anti-self-dual connections (they only differ by the choice 
of orientation on $X$), one can assume without loss of generality that 
$c>0$. Furthermore, a framing for an instanton on $X$ at a point
$p\in X$ is the choice of an isomorphism $E_p\simeq\cpx^n$, where
$E_p$ denotes the fiber of $E$ at $p$.

In the celebrated paper \cite{ADHM}, Atiyah, Drinfeld, Hitchin and Manin
constructed a class of ASD connections on the simplest compact
four dimensional manifold, namely the four dimensional sphere $S^4$.
They also proved that their construction is complete in the sense that
any ASD connection on $S^4$ is gauge equivalent to a connection
constructed by them from a certain algebraic data that depends on the
rank $n$ of the vector bundle $E$ and the instanton number $c$.

More precisely, let $V$ and $W$ be Hermitian vector spaces of dimensions
$c$ and $n$, respectively. The algebraic data of Atiyah, Drinfeld, Hitchin 
and Manin consists of four linear operators:
\begin{equation} \label{adhm.ops}
B_1, B_2 \in {\rm End}(V), \ \ \ \ \ \
i \in {\rm Hom}(W,V), \ \ \ \ \ \ j \in {\rm Hom}(V,W)
\end{equation}
satisfying the following linear relations ($\dagger$ denotes 
Hermitian conjugation):
\begin{eqnarray} 
\label{adhm.eqn1} [ B_1 , B_2 ] + ij & = & 0  \\
\label{adhm.eqn2} [ B_1 , B_1^\dagger ] + [ B_2 , B_2^\dagger ] + ii^\dagger -
j^\dagger j & = & 0
\end{eqnarray}
plus a regularity condition which we describe below.

The ADHM data admits a natural action of the unitary group $U(V)$:
\begin{equation} \label{action}
g(B_1,B_2,i,j) = (gB_1g^{-1},gB_2g^{-1},gi,jg^{-1}), \ \ \ g\in U(V)
\end{equation}
We say that the ADHM datum $(B_1,B_2,i,j)$ is {\em regular} if its
stabilizer subgroup is trivial. Equivalently, $(B_1,B_2,i,j)$ is regular
if and only if it satisfies the following two conditions:

\begin{itemize}
\item {\em stability}: there is no proper subspace $S\subset V$ 
such that $B_k(S)\subset S$ ($k=1,2$) and $i(W)\subset S$;
\item {\em costability}: there is no proper subspace $S\subset V$
such that $B_k(S)\subset S$ ($k=1,2$) and $S\subset \ker j$.
\end{itemize}

The following lemma is probably well known to the experts:

\begin{lemma} \label{st=cost}
Suppose that $(B_1,B_2,i,j)$ satisfies the ADHM equations
(\ref{adhm.eqn1}) and (\ref{adhm.eqn2}). Then $(B_1,B_2,i,j)$
is stable if and only if it is also costable.
\end{lemma}
\begin{proof}
If $(B_1,B_2,i,j)$ is not stable, then by duality on $V$ there
is a proper subspace $S^\perp\subset V$ such that
$B_k^\dagger(S^\perp)\subset S^\perp$ and $S^\perp\subset\ker i^\dagger$.
So restricting (\ref{adhm.eqn2}) to $S^\perp$ and taking the trace, 
we conclude that ${\rm Tr}(j^\dagger j|_{S^\perp}) = 0$. Hence
$S^\perp \subset \ker j$, and $(B_1,B_2,i,j)$ is not costable.

Conversely, $(B_1,B_2,i,j)$ is not costable, take $S^\perp\subset V$ 
nonempty such that $B_k^\dagger(S^\perp)\subset S^\perp$ and 
$S^\perp\subset \ker j$. Restricting (\ref{adhm.eqn2}) to $S^\perp$ 
and taking the trace we conclude that ${\rm Tr}(i i^\dagger|_{S^\perp}) = 0$.
Hence $S^\perp \subset \ker i^\dagger$, and dualizing equation
(\ref{adhm.eqn1}) we see that $[B_1^\dagger,B_2^\dagger]|_{S^\perp}=0$.
Thus $S$ is a proper subspace of $V$ such that
$B_k(S)\subset S$ and $i(W)\subset S$, contradicting stability.
\hfill $\Box$ \end{proof}

We denote the space of regular orbits by:
\begin{equation} \label{moduli}
{\cal M}^{\rm reg}(n,c) = \{ (B_1,B_2,i,j) \ | \
(\ref{adhm.eqn1}) , (\ref{adhm.eqn2}) \} / U(V)
\end{equation}
The main result of
\cite{ADHM} is the following:

\begin{theorem}
The space ${\cal M}^{\rm reg}(n,c)$ is the moduli space of framed 
instantons of rank $n$ and charge $c$ on $S^4$. In other words, 
there is a bijection between points of ${\cal M}^{\rm reg}(n,c)$ 
and gauge equivalence classes of framed ASD connections of rank 
$n$ and charge $c$.
\end{theorem}

%--------------------------------------------------------------------

\subsection{Quantum Haar measure and duality}

For later reference, we now explain the notion of duality on 
free $\cm^{\rm I}_q$- and $\cm^{\rm J}_q$-modules.
We also now
specialize the formal parameter $q$ to a positive real number
for the remainder of this section.

Recall that a Haar functional on a Hopf algebra $\cal A$ is a linear 
functional $H:{\cal A}\seta\cpx$ satisfying the following conditions:
\begin{itemize}
\item bi-invariance: 
$(H\otimes\id)\circ\Delta(a) = (\id\otimes H)\circ\Delta(a) = H(a)$;
\item antipode invariance: $H(\gamma(a))=H(a)$;
\item normalization: $H(1)=1$.
\end{itemize}

\begin{theorem}
\label{thm-h}
There is a unique Haar functional $H$ on the quantum group $SL(2)_q$,
which induces a positive definite Hermitian form on $SU(2)_q$,
namely:
\begin{equation} \label{h}
( g_1 , g_2 ) = H(f^\dagger g), \ \ \ g_1,g_2\in SL(2)_q
\end{equation} \end{theorem}
\begin{proof}
The existence and uniqueness of the Haar functional was established in
\cite{K,MU,W0}. The verification of the fact that $H$ induces a Hermitian
form also can be found in these papers.
\hfill $\Box$ \end{proof}

Now let $U$ be a finite dimensional Hermitian vector space, and let
$(\cdot,\cdot)$ denote
its Hermitian inner product, which is chosen
to be conjugate linear in the first argument. We define the
pairing:
\begin{eqnarray}
(\cdot,\cdot) : 
U\otimes\cm^{\rm I}_q \times U\otimes\cm^{\rm I}_q & \seta & 
\cm^{\rm I}_q \nonumber \\
( v_1\otimes f_1,v_2\otimes f_2 ) & = & 
(v_1,v_2) f_1^\dagger f_2
\label{innerprod}
\end{eqnarray}
From a geometrical point of view, the pairing above plays
the role of a Hermitian metric on a (trivial) vector bundle
over the Euclidean $\real^4$: the pairing of two sections of
the bundle gives a function on $\real^4$.

\begin{proposition} \label{n-deg2}
The pairing (\ref{innerprod}) is non-degenerate, i.e.:
$$ ( \sigma , \sigma ) = 0 \ \ \Leftrightarrow \ \ \sigma=0 $$
\end{proposition}
\begin{proof}
We extend the Haar functional to $\widetilde{SL}(2)_q$ as the homomorphism 
$\tilde{H}:\widetilde{SL}(2)_q\to \cpx[\delta,\delta^{-1}]$ such that 
($n\in\zed$):
$$ \tilde{H} (\delta^n g) =  \delta^n H(g) $$

Now let $\{v_k\}$ be an orthonormal basis for $U$, and 
take $\sigma=\sum_k v_k\otimes f_k$, where $f_k\in\cm^{\rm I}_q$;
we can assume that $f_k=\sum_{\alpha\geq
0}\delta^\alpha f_{k\alpha}$, with
$f_{k\alpha}\in SL_q(2)$.

Then $0=\tilde{H}((\sigma,\sigma))$ forces, by the non degeneracy of
$H$, $f_{k0}=0$, {\it i.e.} $\sigma=\sigma_1\delta$, for some $\sigma_1$.
Then $0=(\sigma,\sigma)=\delta(\sigma_1,\sigma_1)\delta$, and
$(\sigma_1,\sigma_1)=0$ by the invertibility of $\delta$
in $\widetilde{SL_q(2)}$. By iterating this procedure we get that
$\sigma = 0$.

\hfill $\Box$ \end{proof}

Now consider $F_l=U_l\otimes\cm^{\rm I}_q$ as right $\cm^{\rm I}_q$-modules,
where $U_l$ are Hermitian vector spaces, $l=1,2$. Let $F_l^\dagger$ denote
the set of all maps $\mu:F_l\to\cm^{\rm I}_q$ such that
$\mu(\sigma x)=\mu(\sigma)x$, for all $\sigma\in F_l$ and
$x\in\cm^{\rm I}_q$, with the structure of a right $\cm^{\rm I}_q$-module
defined by:
$$ (\mu x) (\sigma) = x^\dagger\mu(\sigma)
\ \ \ \forall \sigma\in F_l $$ 

By Proposition \ref{n-deg2}, the map:
\begin{eqnarray}
{\cal I}_l: F_l & \longrightarrow & F_l^\dagger \label{idn} \\
\sigma & \mapsto & ( \sigma , \cdot ) \nonumber    
\end{eqnarray}
is injective, and it is easy to see that
$\cal I$ must also be surjective. Moreover, notice that
$$ {\cal I}(\sigma x) = x^\dagger {\cal I}(\sigma) $$
In other words, the map $\cal I$ provides an identification
between a free $\cm^{\rm I}_q$-module and its dual.
In particular, one can also regard $\dagger$ as an
involution on $F_l$, defined as follows:
\begin{equation} \label{modinv}
(v\otimes f)^\dagger = \overline{v}\otimes f^\dagger
\end{equation}
and it is easy to see that
$(\sigma x)^\dagger = x^\dagger\sigma^\dagger$ (with this equality
being understood in terms of the bimodule structure of the free
modules $F_l$).

Any map ${\cal L}:F_1\to F_2$ satisfying
${\cal L}(\sigma x)={\cal L}(\sigma)x$ induces a dual map 
${\cal L}^\vee:F_2^\dagger \to F_1^\dagger$ in the usual way:
$$ {\cal L}^\vee (\varphi) = \varphi\circ {\cal L} $$
Finally, the identification (\ref{idn}) yields a map
${\cal L}^\dagger=
{\cal I}_1^{-1}{\cal L}^\vee{\cal I}_2:F_2\to F_1$
with the property:
$$ ( {\cal L}^\dagger \sigma_1 , \sigma_2 ) = 
( \sigma_1 , {\cal L} \sigma_2 ) \ \ \ \ \sigma_l\in F_l $$
In particular, let us consider ${\cal L}=L\otimes x$, where
$L\in{\rm Hom}(U_1,U_2)$ and $x$ means multiplication by
$x\in\cm^{\rm I}_q$ in the left. Then the definition
(\ref{innerprod}) immediately implies that 
${\cal L}^\dagger = L^\dagger\otimes {x^\dagger}$.

\begin{proposition} \label{decomposition}
If ${\cal L}:F_1\to F_2$ is injective, then
$$ F_2 = {\rm Im}{\cal L} \oplus \ker {\cal L}^\dagger $$
\end{proposition}
\begin{proof}
Since ${\cal L}$ is injective, it is easy to see that ${\rm Im}{\cal L}$
is a free submodule of $F_2$. Let $N\subset F_2$ be such that
$F_2={\rm Im}{\cal L} \oplus N$. Define: 
$$ ({\rm Im}{\cal L})^0 =
\{ \psi\in F_2^\dagger \ | \ \psi(\nu)=0 \ \forall \nu\in {\rm Im}{\cal L} \} $$
One can then show that
$N^\dagger \simeq ({\rm Im}{\cal L})^0 \simeq \ker{\cal L}^\vee$.
The desired decomposition now follows from the identification 
$N\simeq N^\dagger$. 
\hfill $\Box$ \end{proof}

%--------------------------------------------------------------------

\subsection{Construction of quantum instantons}

First of all, we must explain precisely what we mean by a
quantum instanton on the quantum affine Minkowski space.

\begin{definition} \label{nc.instanton}
A quantum instanton on $S^{\rm I}_q$ is a triple 
$(E,\nabla,\dagger)$ consisting of:
\begin{itemize}
\item a finitely generated, projective right $\cm^{\rm I}_q$-module $E$ 
equipped with an involution $\dagger: E \rightarrow E$
satisfying $(\sigma x)^\dagger = x^\dagger \sigma^\dagger$
for all $x\in\cm^{\rm I}_q$ and $\sigma\in E$;
\item anti-self-dual connection 
$\nabla:E\seta E\otimes\Omega^1_{\cm^{\rm I}_q}$ 
which is compatible with the involution $\dagger$,
i.e. $\nabla \dagger = \dagger \nabla$.
\end{itemize} \end{definition}

Quantum instantons on $S^{\rm J}_q$ are similarly defined. Moreover,
we also define the notion of {\em consistency} between quantum 
instantons on $S^{\rm I}_q$ and $S^{\rm J}_q$:

\begin{definition} \label{consistency}
Quantum instantons $(E_{\rm I},\nabla_{\rm I},\dagger_{\rm I})$ 
on $S^{\rm I}_q$ and $(E_{\rm J},\nabla_{\rm J},\dagger_{\rm J})$
on $S^{\rm J}_q$ are said to be {\em consistent} if:
\begin{itemize}
\item there is an isomorphism
\begin{equation} \label{nc.mod.glue.map}
\Gamma: E_{\rm I} [\delta^{-1}] \seta E_{\rm J}[ \delta ]
\end{equation}
such that $\Gamma(\sigma f)=\Gamma(\sigma)\eta(f)$, for all
$\sigma\in E_{\rm I}[\delta^{-1}]$ and $f\in \cm^{\rm I}_q[\delta^{-1}]$;
\item $\nabla_{\rm J} \Gamma = \Gamma \nabla_{\rm I}$;
\item $\dagger_{\rm J} \Gamma = \Gamma \dagger_{\rm I}$.
\end{itemize} \end{definition}

Recall that $\eta$ is the isomorphism
$\cm^{\rm I}_q[\delta^{-1}] \seta \cm^{\rm J}_q[\delta]$
described in (\ref{nc.glueing}). Geometrically, the consistency
condition means that the quantum instantons $\nabla_{\rm I}$
and $\nabla_{\rm J}$ coincide in the ``intersection'' variety
$\cm^{\rm I}_q[\delta^{-1}]\simeq\cm^{\rm J}_q[\delta]$, up to a gauge
transformation.

The goal of this Section is to convert ADHM data into a consistent pair
of quantum instantons, in close analogy with the classical case.

\paragraph{Quantum Instantons on quantum Minkowski space-time}
As before, let $V$, $W$ denote Hermitian vector spaces of dimension $c$
and $n$, respectively. Let $\tilde{W}=V\oplus V\oplus W$.

Let $(B_1,B_2,i,j)$ be an ADHM datum, as in (\ref{adhm.ops}). 
We start by considering the following sequence of free $\cm_q^{\rm I}$-modules:
$$ V\otimes\cm_q^{\rm I} \stackrel{\alpha_{\rm I}}{\longrightarrow}
\tilde{W}\otimes\cm_q^{\rm I} \stackrel{\beta_{\rm I}}{\longrightarrow}
V\otimes\cm_q^{\rm I} $$
where the maps $\alpha_{\rm I}$ and $\beta_{\rm I}$ are given by:
\begin{equation}\label{alpha}
\alpha_{\rm I} = \left( \begin{array}{c}
B_1\otimes\id - \id\otimes x_{21'} \\
B_2\otimes\id - \id\otimes x_{22'} \\
j\otimes\id
\end{array} \right) \end{equation}
and
\begin{equation}\label{beta}
\beta_{\rm I} = \left( \begin{array}{ccc}
-B_2\otimes\id + \id\otimes x_{22'} \ \ &
B_1\otimes\id - \id\otimes x_{21'} \ \ &
i\otimes\id
\end{array} \right) \end{equation}

The induced dual maps
$\beta_{\rm I}^\dagger:V\otimes\cm_q^{\rm I}\to\tilde{W}\otimes\cm_q^{\rm I}$ and 
$\alpha_{\rm I}^\dagger:\tilde{W}\otimes\cm_q^{\rm I}\to V\otimes\cm_q^{\rm I}$ 
are given by:
\begin{equation}\label{betadual}
\beta_{\rm I}^\dagger = \left( \begin{array}{c}
-B_2^\dagger\otimes\id + \id\otimes x_{11'} \\
B_1^\dagger\otimes\id + \id\otimes  x_{12'} \\
i^\dagger\otimes\id
\end{array} \right) \end{equation}
and
\begin{equation}\label{alphadual}
\alpha_{\rm I}^\dagger = \left( \begin{array}{lcr}
B_1^\dagger\otimes\id + \id\otimes x_{12'} \ \ &
B_2^\dagger\otimes\id - \id\otimes x_{11'} \ \ &
j^\dagger\otimes\id
\end{array} \right) \end{equation}

\begin{proposition} \label{compo}
\begin{enumerate}
\item $\beta_{\rm I}\alpha_{\rm I}=0$ if and only if $[B_1,B_2]+ij=0$.
\item $\beta_{\rm I}\beta_{\rm I}^\dagger=
\alpha_{\rm I}^\dagger\alpha_{\rm I}$
if and only if 
$[B_1,B_1^\dagger]+[B_2,B_2^\dagger]+ii^\dagger-j^\dagger j=0$.
\end{enumerate}
\end{proposition}
\begin{proof}
It is easy to check that:
$$ \beta_{\rm I}\alpha_{\rm I} = 
\left( [B_1,B_2]+ij \right) \otimes\id + 
\id\otimes [x_{22'},x_{21'}] $$
Since $[x_{22'},x_{21'}]=0$ (see equation (\ref{xcommut1})),
the first statement follows.

Another straightforward calculation reveals that:
$$ \beta_{\rm I}\beta_{\rm I}^\dagger -\alpha_{\rm I}^\dagger\alpha_{\rm I} =
\left( [B_1,B_1^\dagger]+[B_2,B_2^\dagger]+ii^\dagger-j^\dagger j \right) \otimes \id
+  \id \otimes ([x_{11'},x_{22'}]-[x_{12'},x_{21'}]) $$
and the second statement follows from equation (\ref{xcommut3}). 
\hfill $\Box$ \end{proof} 

\begin{proposition} \label{inj-surj}
\begin{enumerate}
\item $\alpha_{\rm I}$ is injective.
\item $\beta_{\rm I}$ is surjective if and only if 
$(B_1,B_2,i,j)$ is stable.
\end{enumerate}
\end{proposition}
\begin{proof}
Let ${\cal X}_1=\cpx[x_{11'},x_{12'}]$ and ${\cal X}_2=\cpx[x_{21'},x_{22'}]$.
Then $\cm^{\rm I}_q$ can be represented
in the following way:
$$ \cm^{\rm I}_q = {\cal X}_1 \otimes {\cal X}_2 / 
(\ref{xcommut3}),(\ref{xcommut2}) $$
Clearly, $\alpha_{\rm I}$ and $\beta_{\rm I}$ can be restricted to maps
$$ V\otimes{\cal X}_2 \stackrel{\alpha_{\rm I}}{\longrightarrow}
\tilde{W}\otimes{\cal X}_2 \stackrel{\beta_{\rm I}}{\longrightarrow}
V\otimes{\cal X}_2
$$
It then follows from Nakajima (see \cite{N2}, lemma 2.7) that
$\alpha_{\rm I}|_{{\cal X}_2}$ is injective and that $\beta_{\rm I}|_{{\cal X}_2}$
is surjective if and only if $(B_1,B_2,i,j)$ is stable.
Since $\alpha_{\rm I}=\id_{{\cal X}_1}\otimes\alpha_{\rm I}|_{{\cal X}_2}$ and 
$\beta_{\rm I}=\id_{{\cal X}_1}\otimes\beta_{\rm I}|_{{\cal X}_2}$,
we conclude that $\alpha_{\rm I}$ is also injective, and that 
$\beta_{\rm I}$ is surjective if and only if $(B_1,B_2,i,j)$ is stable.
\hfill $\Box$ \end{proof}

Thus we have the short sequence of free $\cm^{\rm I}_q$-modules:
\begin{equation} \label{nc.monad2}
0 \seta V\otimes\cm_q^{\rm I} \stackrel{\alpha_{\rm I}}{\longrightarrow}
\tilde{W}\otimes\cm_q^{\rm I} \stackrel{\beta_{\rm I}}{\longrightarrow}
V\otimes\cm_q^{\rm I} \seta 0
\end{equation}
which is exact on the first and last terms. Its middle cohomology 
$E_{\rm I}={\rm ker}\beta_{\rm I}/{\rm Im}\alpha_{\rm I}$ is then a well
defined right $\cm^{\rm I}_q$-module since both $\beta_{\rm I}$ and 
$\alpha_{\rm I}$ are right linear. Clearly $E_{\rm I}$ is torsion-free,
since it is a submodule of a free module; we argue that it is projective. 

Furthermore, notice that the pairing (\ref{innerprod}) also induces a
pairing $E_{\rm I}\times E_{\rm I}\to \cm^{\rm I}_q$.

\begin{proposition} \label{inj-surj2}
\begin{enumerate}
\item $\beta_{\rm I}^\dagger$ is injective.
\item $\alpha_{\rm I}^\dagger$ is surjective if and only if 
$(B_1,B_2,i,j)$ is costable.
\end{enumerate}
\end{proposition}
\begin{proof}
Same argument as Proposition \ref{inj-surj}, just replacing ${\cal X}_2$
by ${\cal X}_1$ and vice-versa.
\hfill $\Box$ \end{proof}

\begin{proposition} \label{iso}
$\xi_{\rm I}=\beta_{\rm I}\beta_{\rm I}^\dagger=\alpha_{\rm I}^\dagger\alpha_{\rm I}$ 
is an isomorphism if and only if $(B_1,B_2,i,j)$ is stable.
\end{proposition}

\begin{proof}
Applying Proposition \ref{decomposition} to the map $\beta_{\rm I}^\dagger$, 
we get that:
\begin{equation} \label{decomp-beta}
\tilde{W}\otimes\cm_q^{\rm I} = 
{\rm Im}\beta_{\rm I}^\dagger \oplus \ker\beta_{\rm I}
\end{equation} 
In particular, ${\rm Im}\beta_{\rm I}^\dagger \cap \ker \beta_{\rm I} =\{0\}$
so that the injectivity of $\xi_{\rm I}$ follows from the
injectivity of $\beta_{\rm I}^\dagger$.

Now if $(B_1,B_2,i,j)$ is stable then $\beta_{\rm I}$ is surjective, so given
$\nu\in V\otimes\cm_q^{\rm I}$, there is $\mu\in\tilde{W}\otimes\cm_q^{\rm I}$
such that $\beta_{\rm I}(\mu)=\nu$. By the decomposition
(\ref{decomp-beta}), we know that $\mu=\mu'+\mu''$ for some 
$\mu'\in {\rm Im}\beta_{\rm I}^\dagger$ and
$\mu''\in \ker\beta_{\rm I}$. Let $\mu'=\beta_{\rm I}^\dagger(\nu')$.
Then $\nu=\beta_{\rm I}(\mu)=\beta_{\rm I}(\mu')=
\beta_{\rm I}(\beta_{\rm I}^\dagger(\nu'))$, hence $\xi_{\rm I}$ is also
surjective.

Conversely, if $(B_1,B_2,i,j)$ is not stable, then $\beta_{\rm I}$ is not
surjective, so $\beta_{\rm I}\beta_{\rm I}^\dagger$ is not surjective either,
and $\xi_{\rm I}$ is not an isomorphism.
\hfill $\Box$ \end{proof}

Now consider the {\em Dirac operator}:  
\begin{eqnarray} 
{\cal D}_{\rm I} :
\tilde{W}\otimes\cm_q^{\rm I} & \seta & (V\oplus V)\otimes\cm_q^{\rm I} \nonumber \\ 
\label{dirac} {\cal D}_{\rm I} & = & 
\left( \begin{array}{c} \beta_{\rm I} \\ \alpha_{\rm I}^\dagger \end{array} \right) 
\end{eqnarray} 
Moreover, we define the {\em Laplacian}:
$$ \Xi_{\rm I} :  (V\oplus V)\otimes\cm_q^{\rm I} \seta 
(V\oplus V)\otimes\cm_q^{\rm I} $$ 
\begin{equation} \label{nc.laplacian}
\Xi_{\rm I} = {\cal D}_{\rm I}{\cal D}_{\rm I}^\dagger = 
\left( \begin{array}{cc} \beta_{\rm I}\beta_{\rm I}^\dagger & 0 \\ 
0 & \alpha_{\rm I}^\dagger\alpha_{\rm I} \end{array} \right)
\end{equation} 
Proposition \ref{iso} implies that $\Xi_{\rm I}=\xi_{\rm I}\id_{V\oplus V}$ 
and $\Xi_{\rm I}$ is an isomorphism too.  We can then define the projection map:  
$$ P_{\rm I} : \tilde{W}\otimes\cm_q^{\rm I} \seta E_{\rm I} $$ 
\begin{equation} \label{proj}
P_{\rm I} = \id - {\cal D}_{\rm I}^\dagger\Xi_{\rm I}^{-1}{\cal D}_{\rm I}
\end{equation} 
For the next two Propositions, we assume that $(B_1,B_2,i,j)$ is stable.  We have:

\begin{proposition} \label{proj2}
$E_{\rm I} \simeq {\rm ker}{\cal D}_{\rm I}$.  In particular,
$E_{\rm I}$ is a projective right $\cm_q^{\rm I}$-module.  
\end{proposition}
\begin{proof} 
Given $\psi\in{\rm ker}\beta_{\rm I}$, we show that there is a unique
$\nu\in V\otimes\cm_q^{\rm I}$ such that 
$\psi'=\psi+\alpha_{\rm I}(\nu)\in\ker{\cal D}_{\rm I}$, i.e.  
$\beta_{\rm I}(\psi')=\alpha_{\rm I}^\dagger(\psi')=0$.  Indeed:  
$$ \beta_{\rm I}(\psi')=\beta_{\rm I}\alpha_{\rm I}(\nu)=0 $$ 
$$ \alpha_{\rm I}^\dagger(\psi')=0 \Leftrightarrow 
\alpha_{\rm I}^\dagger(\psi)= -\alpha_{\rm I}^\dagger\alpha_{{\rm I}}(\nu) $$ 
but $\xi_{\rm I}=\alpha_{\rm I}^\dagger\alpha_{\rm I}$ is an isomorphism,
thus $\nu=\xi_{\rm I}^{-1}\alpha_{\rm I}^\dagger(\psi)$, as desired.

Finally, it is easy to see that
given $\nu\in\tilde{W}\otimes\cm_q^{\rm I}$,
there are unique $\psi\in\ker{\cal D}_{\rm I}$ and 
$\varphi\in{\rm Im}{\cal D}_{\rm I}^\dagger$ such that
$\nu=\psi+\varphi$. Indeed, just take $\psi=P_{\rm I}\nu$
and $\varphi={\cal D}_{\rm I}^\dagger\Xi_{\rm I}^{-1}{\cal D}_{\rm I}\nu$.

In other words, we conclude that
$E_{\rm I}\oplus{\rm Im}{\cal D}_{\rm I}^\dagger=
\tilde{W}\otimes\cm_q^{\rm I}$, which implies that
$E_{\rm I}$ is projective as a $\cm_q^{\rm I}$-module. 
\hfill $\Box$ \end{proof}

Note also that $E_{\rm I}$ is finitely generated and has rank $n=\dim W$. 

To define the connection, let 
$\iota_{\rm I}:E_{\rm I}\seta\tilde{W}\otimes\cm_q^{\rm I}$ 
denote the natural inclusion and $d:\cm_q^{\rm I}\seta\Omega^1_{\cm_q^{\rm I}}$
denote the quantum de Rham operator. We define $\nabla_{\rm I}$ 
via the composition:
$$ \xymatrix{
E_{\rm I} \ar[r]^{\iota_{\rm I}} &
\tilde{W}\otimes\cm_q^{\rm I} \ar[r]^{\id\otimes d} &
\tilde{W}\otimes\Omega^1_{\cm_q^{\rm I}} \ar[r]^{P_{\rm I}\otimes\id} &
E_{\rm I}\otimes_{\cm_q^{\rm I}}\Omega^1_{\cm_q^{\rm I}}
} $$

\begin{proposition} \label{asd}
$F_{\nabla_{\rm I}}$ is anti-self-dual.
\end{proposition}
\begin{proof}
Note that $F_{\nabla_{\rm I}}=\nabla_{\rm I}\nabla_{\rm I}=P_{\rm I}dP_{\rm I}d$;
%is given by the composition:
%$$ \xymatrix{
%E_{\rm I} \ar[r]^{\iota_{\rm I}} &
%\tilde{W}\otimes\cm_q^{\rm I} \ar[r]^{\id\otimes d} &
%\tilde{W}\otimes\Omega^1_{\cm_q^{\rm I}} \ar[r]^{P_{\rm I}\otimes\id} &
%E_{\rm I}\otimes_{\cm_q^{\rm I}}\Omega^1_{\cm_q^{\rm I}} \ar[r]^{\id\otimes d} & \\
%& \ar[r]^{\id\otimes d} & \tilde{W}\otimes\Omega^2_{\cm_q^{\rm I}} \ar[r]^{P_{\rm I}\otimes\id} &
%E_{\rm I}\otimes_{\cm_q^{\rm I}}\Omega^2_{\cm_q^{\rm I}}
%} $$
therefore given $e\in E_{\rm I}$ we have:
\begin{eqnarray*}
F_{\nabla_{\rm I}} e & = &
P_{\rm I} \left( d ( \id - {\cal D}_{\rm I}^\dagger\Xi_{\rm I}^{-1}{\cal D}_{\rm I} ) de \right) = 
P_{\rm I} \left( d {\cal D}_{\rm I}^\dagger\Xi_{\rm I}^{-1} (d{\cal D}_{\rm I}) e \right) = \\
& = & P_{\rm I} \left( (d{\cal D}_{\rm I}^\dagger)\Xi_{\rm I}^{-1} (d{\cal D}_{\rm I}) e +
{\cal D}_{\rm I}^\dagger d(\Xi_{\rm I}^{-1} (d{\cal D}_{\rm I}) e) \right) = \\
& = & P_{\rm I} \left( (d{\cal D}_{\rm I}^\dagger)\Xi_{\rm I}^{-1} (d{\cal D}_{\rm I}) e \right)
\end{eqnarray*}
for $P_{\rm I} \left({\cal D}_{\rm I}^\dagger d(\Xi_{\rm I}^{-1} (d{\cal D}_{\rm I}) e) \right) = 0$.
Since $\Xi_{\rm I}^{-1}=\xi_{\rm I}^{-1}\id$, we conclude that  
$F_{\nabla_{\rm I}}$ is proportional to
$d{\cal D}_{\rm I}^\dagger\wedge d{\cal D}_{\rm I}$, as a 2-form.
 
It is then a straightforward calculation to show that each entry of 
$d{\cal D}_{\rm I}^\dagger\wedge d{\cal D}_{\rm I}$ belongs to
$\Omega^{2,-}_{\cm_q^{\rm I}}$; indeed:
$$ d{\cal D}_{\rm I}^\dagger\wedge d{\cal D}_{\rm I} = 
\left( \begin{array}{lr} 
dx_{11'} \ \ & \ \ -dx_{21'} \\
dx_{12'} \ \ & \ \ -dx_{22'} \\ 0 & 0
\end{array} \right) \wedge
\left( \begin{array}{lrc}
dx_{22'} \ \ & \ \  -dx_{21'} & \ 0 \\
dx_{12'} \ \ & \ \ -dx_{11'} & \ 0
\end{array} \right) = $$
$$ = \left( \begin{array}{ccc}
dx_{11'}dx_{22'}-dx_{21'}dx_{12'} \ \ &
\ \ -dx_{11'}dx_{21'}+dx_{21'}dx_{11'} & \ 0 \\
dx_{12'}dx_{22'}-dx_{22'}dx_{12'} \ \ &
\ \ -dx_{12'}dx_{21'}+dx_{22'}dx_{11'} & \ 0 \\
0 & 0 & \ 0 
\end{array} \right)$$
Applying the commutation relations (\ref{xnc2forms}), we obtain:
$$ d{\cal D}_{\rm I}^\dagger\wedge d{\cal D}_{\rm I} = 
\left( \begin{array}{ccc}
dx_{11'}dx_{22'}+dx_{12'}dx_{21'} &
-2dx_{11'}dx_{21'} & \ 0 \\
2dx_{12'}dx_{22'} & 
-(dx_{11'}dx_{22'}+dx_{12'}dx_{21'})& \ 0 \\
0 & 0 & \ 0 
\end{array} \right)$$
Comparing with (\ref{nc-asd2f}), we have proved our claim.
\hfill $\Box$ \end{proof}

\paragraph{Gauge equivalence}
We show that if $(B_1,B_2,i,j)$ and $(B_1',B_2',i',j')$ are
equivalent ADHM data, then the respective pairs $(E_{\rm I},\nabla_{\rm I})$
and $(E_{\rm I}',\nabla'_{\rm I})$ are gauge equivalent, in the sense
that there is a $\cm^{\rm I}_q$-isomorphism $G:E_{\rm I}'\seta E_{\rm I}$
such that $\nabla_{\rm I}'=G^{-1}\nabla_{\rm I}G$.

To do that, recall that $(B_1,B_2,i,j)$ and $(B_1',B_2',i',j')$ are
equivalent if there exists $g\in{\rm U}(V)$ such that:
$$ \begin{array}{ccccc}
B_k' = g B_k g^{-1},\ k=1,2 & \ \ \ &
i' = g i & \ \ \ & j' = j g^{-1}
\end{array} $$
Let $G\in U(\tilde{W})$ be given by $g\times g\times\id_{W}$.
It is then easy to check that the following diagram is commutative:
\begin{equation} \label{geqv}
\xymatrix{
0 \ar[r] & V\otimes\cm^{\rm I}_q \ar[r]^{\alpha_{\rm I}'}
\ar[d]^{g\otimes\id} &
\tilde{W}\otimes\cm^{\rm I}_q \ar[r]^{\beta_{\rm I}'}
\ar[d]^{G\otimes\id} &
V\otimes\cm^{\rm I}_q \ar[r] \ar[d]^{g\otimes\id} & 0 \\
0 \ar[r] & V\otimes\cm^{\rm I}_q \ar[r]^{\alpha_{\rm I}} &
\tilde{W}\otimes\cm^{\rm I}_q \ar[r]^{\beta_{\rm I}}&
V\otimes\cm^{\rm I}_q \ar[r] & 0
} \end{equation}
Therefore the modules 
$E_{\rm I}={\rm ker}\beta_{\rm I}/{\rm Im}\alpha_{\rm I}$ and
$E_{\rm I}'={\rm ker}\beta_{\rm I}'/{\rm Im}\alpha_{\rm I}'$ are
isomorphic; indeed, it is easy to see that $G$ maps $E_{\rm I}'$
onto $E_{\rm I}$ (regarded as submodules of 
$\tilde{W}\otimes\cm^{\rm I}$). We shall also denote by
$G$ the induced isomorphism $E_{\rm I}'\seta E_{\rm I}$).

Now denote by $\iota_{\rm I}':E_{\rm I}'\seta\tilde{W}\otimes\cm^{\rm I}_q$
the
inclusion and by 
$P_{\rm I}':\tilde{W}\otimes\cm^{\rm I}_q\seta E_{\rm I}'$
the projection (\ref{proj}). Clearly, $\iota_{\rm I}'=G^{-1}\iota_{\rm I} G$
and
$P_{\rm I}'=G^{-1} P_{\rm I} G$. In addition, we have:
\begin{eqnarray*}
\nabla_{\rm I}' & = & P_{\rm I}' d \iota_{\rm I}' = 
G^{-1} P_{\rm I} G d G^{-1}\iota_{\rm I} G = \\
& = & G^{-1} P_{\rm I} \left( GdG^{-1}(\iota_{\rm I} G) + 
d \iota_{\rm I} G \right) = \\
& = & G^{-1} P_{\rm I} d \iota_{\rm I} G = G^{-1}\nabla_{\rm I}G 
\end{eqnarray*}
since $G$ acts as the identity on $\cm^{\rm I}_q$, so that $dG^{-1}=0$.

\paragraph{Real structure}
Dualizing the monad (\ref{nc.monad2}) and using the identifications
$(V\otimes\cm_q^{\rm I})^\dagger\simeq V\otimes\cm_q^{\rm I}$ and
$(\tilde{W}\otimes\cm_q^{\rm I})^\dagger\simeq\tilde{W}\otimes\cm_q^{\rm I}$
one obtains:
$$ 0 \seta V\otimes\cm_q^{\rm I} \stackrel{\beta_{\rm I}^\dagger}{\seta}
\tilde{W}\otimes\cm_q^{\rm I} 
\stackrel{\alpha_{\rm I}^\dagger}{\seta}
V\otimes\cm_q^{\rm I} \seta 0 $$
From Proposition \ref{inj-surj2},
this monad is again exact at the first and last terms.
Let us denote its cohomology by $E_{\rm I}^\dagger$, which also has the 
structure of a right $\cm_q^{\rm I}$-module. Moreover, via the procedure
in the proof of Proposition \ref{proj2}, $E_{\rm I}^\dagger$
can be
identified with the kernel of the map:
$$ \left( \begin{array}{c} 
\alpha_{\rm I}^\dagger \\ \beta_{\rm I}
\end{array} \right) : \tilde{W}\otimes\cm_q^{\rm I}
\seta (V\oplus V)\otimes\cm_q^{\rm I} $$
which is clearly isomorphic to ${\rm ker}{\cal D}_{\rm I}\simeq E_{\rm I}$. 
Therefore, the involution 
$\dagger: \tilde{W}\otimes\cm_q^{\rm I}\to\tilde{W}\otimes\cm_q^{\rm I}$
induces a map $\dagger:E_{\rm I}\seta E_{\rm I}$; the desired property
follows easily from (\ref{modinv}).
To check the compatibility of $\dagger$
with the connection $\nabla_{\rm I}$,
note that:
$$ \nabla_{\rm I}(e^\dagger) = P_{\rm I}d(e^\dagger) =
P_{\rm I}(de)^\dagger = (P_{\rm I}de)^\dagger = (\nabla_{\rm I}e)^\dagger $$

Summing up the work done so far, we have constructed a well-defined map
from the set of equivalence classes of {\em classical} ADHM data to the
set of gauge equivalence classes of {\em quantum} instantons on quantum
Minkowski space-time $\cm_q^{\rm I}$, in the sense of Definition
\ref{nc.instanton}.

\paragraph{Connection matrix}
Finally, let us describe the connection matrix associated with the connection
$\nabla_{\rm I}$ given above. To do that, let $\{\sigma_k,\rho^k\}_{k=1}^{r}$ 
be a dual basis for $E_{\rm I}$. Let also $\{w_k\}_{k=1}^{r}$ be an orthonormal 
basis for $W$. These choices induce a natural map:
\begin{eqnarray*}
\Psi : W\otimes\cm^{\rm I}_q & \seta & 
\tilde{W}\otimes\cm^{\rm I}_q \\
\Psi(w_k\otimes f) & = & \iota_{\rm I}(\sigma_k)f
\end{eqnarray*}
extended by linearity. Then $\Psi$ has the property ${\cal D}_{\rm I}\circ\Psi=0$:
$$ {\cal D}_{\rm I}\circ\Psi (\sum_k w_k\otimes f_k) =
{\cal D}_{\rm I} (\sum_k f_k  \iota(\sigma_k)) =
\sum f_k  {\cal D}_{\rm I}(\iota(\sigma_k)) = 0 $$
The map $\Psi$ is clearly injective, so that 
$\Psi^\dagger\Psi:W\otimes\cm^{\rm I}_q\seta W\otimes\cm^{\rm I}_q$
is an isomorphism. Moreover, the basis $\{w_k\}_{k=1}^{r}$ can be chosen
such that
$\Psi^\dagger\Psi=\id$. 

On the other hand, consider the map:
\begin{eqnarray*}
\rho : E_{\rm I} & \seta & 
W\otimes\cm^{\rm I}_q \\
\rho(e) & = & \sum_k w_k\otimes \rho^k(e)
\end{eqnarray*}
Then $\rho P_{\rm I} \Psi = \id$, thus $\rho P_{\rm I} = \Psi^\dagger$.

Recalling that the entries of the connection matrix are given by
$A_k^l = \rho^l\otimes\id (\nabla_{\rm I} \sigma_k)$, we have:
$$ A_k^l = \rho^l \left( P _{\rm I} d\Psi (w_k\otimes 1) \right) $$
But the right hand side are just the matrix coefficients of
$\Psi^\dagger d\Psi$, so that $A=\Psi^\dagger d\Psi$, as in the classical
ADHM construction.

\paragraph{Quantum instantons on $S_q^{\rm J}$}
Consider the monad
\begin{equation} \label{nc.monad3}
0 \seta V\otimes\cm_q^{\rm J} \stackrel{\alpha_{\rm J}}{\seta}
\tilde{W}\otimes\cm_q^{\rm J} \stackrel{\beta_{\rm J}}{\seta}
V\otimes\cm_q^{\rm J} \seta 0
\end{equation}
with the maps:
$$ \alpha_{\rm J} = \left( \begin{array}{c}
B_1\otimes\id - \id\otimes  y_{12'} \\
B_2\otimes\id - \id\otimes  y_{22'} \\
j\otimes\id
\end{array} \right) $$
and
$$ \beta_{\rm J} = \left( \begin{array}{ccc}
-B_2\otimes\id + \id\otimes  y_{22'} &
 B_1\otimes\id - \id\otimes  y_{12'} &
i\otimes\id
\end{array} \right) $$
so that $\alpha_{\rm J} = (\id\otimes\eta) \alpha_{\rm I} (\id\otimes\eta^{-1})$
and $\beta_{\rm J} = (\id\otimes\eta) \beta_{\rm I} (\id\otimes\eta^{-1})$ 

It is again easy to check that $\beta_{\rm J}\alpha_{\rm J}=0$, that
$\alpha_{\rm J}$ is injective, and that $\beta_{\rm J}$ is surjective
if and only if $(B_1,B_2,i,j)$ is stable. Thus, the cohomology of
(\ref{nc.monad3}), denoted by $E_{\rm J}$,
is a projective right
$\cm_q^{\rm J}$-module. The connection $\nabla_{\rm J}$
and real
structure $\dagger_{\rm J}:E_{\rm J}\seta E_{\rm J}$ can be similarly 
defined, and we obtain a quantum instanton on $\cm^{\rm J}_q$. 

\paragraph{Consistency}
If the modules $E_{\rm I}$ and $E_{\rm J}$ are constructed as above, 
the consistency map (\ref{nc.mod.glue.map}) arises in the following way.
Consider the diagram:
\begin{equation} \label{nc.diag1} \xymatrix{
0 \ar[r] &
V\otimes\cm^{\rm I}_q[\delta^{-1}] \ar[r]^{\alpha_{\rm I}}
\ar[d]^{\id_V\otimes\eta} &
\tilde{W}\otimes\cm^{\rm I}_q[\delta^{-1}] \ar[r]^{\beta_{\rm I}}
\ar[d]^{\id_{\tilde{W}}\otimes\eta} &
V\otimes\cm^{\rm I}_q[\delta^{-1}] \ar[r] \ar[d]^{\id_V\otimes\eta} & 0 \\
0 \ar[r] &
V\otimes\cm^{\rm J}_q[\delta] \ar[r]^{\alpha_{\rm J}} &
\tilde{W}\otimes\cm^{\rm J}_q[\delta] \ar[r]^{\beta_{\rm J}}&
V\otimes\cm^{\rm J}_q[\delta] \ar[r] & 0 } 
\end{equation}
Clearly, the cohomology of the first row is $E_{\rm I}[\delta^{-1}]$,
while the cohomology of the second row is $E_{\rm J}[\delta]$.
Moreover, the diagram is commutative. Therefore, the isomorphism 
$\id_{\tilde{W}}\otimes\eta$ induces an isomorphism 
$\Gamma:E_{\rm I}[\delta^{-1}]\seta E_{\rm J}[\delta]$,
as required in Definition \ref{consistency}.

To establish the consistency between the connections $\nabla_{\rm I}$ and 
$\nabla_{\rm J}$, it is enough to show that the connection matrices 
$A_{\rm I}$ and $A_{\rm J}$ are related via a gauge transformation
on the ``intersection algebra'' 
$\cm^{\rm I}_q[\delta^{-1}]\stackrel{\eta}{\seta}\cm^{\rm J}_q[\delta]$. 
Indeed, fix a trivialization:
$\Psi_{\rm I}$ of $E_{\rm I}$ such that 
$\Psi_{\rm I}^\dagger \Psi_{\rm I}=\id$; consider the diagram:
$$ \xymatrix{
0 \ar[r] &
W\otimes\cm^{\rm I}_q[\delta^{-1}] \ar[r]^{\Psi_{\rm I}}
\ar[d]^{\id_V\otimes\eta} &
\tilde{W}\otimes\cm^{\rm I}_q[\delta^{-1}] \ar[r]^{{\cal D}_{\rm I}}
\ar[d]^{\id_{\tilde{W}}\otimes\eta} &
V\otimes\cm^{\rm I}_q[\delta^{-1}] \ar[r] \ar[d]^{\id_V\otimes\eta} & 0 \\
0 \ar[r] &
V\otimes\cm^{\rm J}_q[\delta] \ar@{-->}[r]&
\tilde{W}\otimes\cm^{\rm J}_q[\delta] \ar[r]^{{\cal D}_{\rm J}}&
V\otimes\cm^{\rm J}_q[\delta] \ar[r] & 0 } $$
The commutativity of the second square on the above follows from the
commutativity of the diagram (\ref{nc.diag1}). This means that 
$\Psi_{\rm J}=(\id_{\tilde{W}}\otimes\eta)\Psi_{\rm I}(\id_V\otimes\eta^{-1})$
is a trivialization of $E_{\rm J}$. To simplify notation, let us simply use 
$\eta$ to denote $\id\otimes\eta$. Hence:
\begin{eqnarray*}
A_{\rm J} & = & 
(\eta\Psi_{\rm I}\eta^{-1})^\dagger d_{\rm J} (\eta\Psi_{\rm I}\eta^{-1}) =
\eta\Psi_{\rm I}^\dagger\eta^{-1} \eta d_{\rm I} (\Psi_{\rm I}\eta^{-1}) = \\
& = & \eta (\Psi_{\rm I}^\dagger d_{\rm I} \Psi_{\rm I})\eta^{-1} +
\eta \Psi_{\rm I}^\dagger \Psi_{\rm I} d_{\rm I} (\eta^{-1}) =
\eta A_{\rm I} \eta^{-1} + \eta d_{\rm I} \eta^{-1}
\end{eqnarray*}
where $d_{\rm I}$ and $d_{\rm J}$ denote the de Rham operators on
$\cm_{\rm I}$ and $\cm_{\rm J}$, respectively. 

In other words, 
$\nabla_{\rm J}=\Gamma\nabla_{\rm I}\Gamma^{-1}$. 
We also used the fact that $(\eta^{-1})^\dagger=\eta$ and that 
$d_{\rm J} \eta = \eta d_{\rm I}$.

Finally, recall that $\dagger\eta=\eta\dagger$. Since $\Gamma$ is
induced from $\eta$ and $\dagger$ is induced from $\dagger$, we conclude
that indeed: $\dagger_{\rm J}\Gamma = \Gamma\dagger_{\rm I}$.

We sum up the work done in this Section in the following statement,
which motivated the title of this paper:

\begin{theorem} \label{oneway}
There exists a well-defined map from the set of equivalence classes of
regular ADHM data to the moduli space of gauge equivalence classes of
consistent
pairs of quantum instantons.
\end{theorem}

%---------------------------------------------------------------
%---------------------------------------------------------------

\section{Quantum Penrose transform and further perspectives} \label{conc}

\subsection{Completeness conjecture and the quantum Penrose transform}

As we mentioned in Introduction, we conjecture that all
anti-self-dual connections on $\cm^{\rm I}_q$ are gauge equivalent 
to the ones produced above. In other words, the map given by Theorem
\ref{oneway} is invertible: given a consistent pair of quantum
instantons, there is an ADHM datum $(B_1,B_2,i,j)$ such that the consistent
pair can be reconstructed from $(B_1,B_2,i,j)$ via the procedure above.

As a consequence of this conjecture, we are able to conclude that the
moduli space of {\em quantum} instantons actually coincides with the
moduli space of {\em classical} instantons, therefore fully
justifying the title of this paper.

The key ingredient in the proof of the classical version of this
conjecture is Penrose's twistor diagram (also termed the 
{\it flat self-duality diagram} by Manin \cite{Ma}):
$$ \xymatrix{
& \mathbf{F}_{1,2}(\st) \ar[dl]^\mu \ar[dr]_\nu & \\
\mathbf{P}(\st) & & \mathbf{G}_2(\st)=\m
} $$
where $\mathbf{F}_{1,2}(\st)$ denotes the flag manifold of lines within planes
in $\st$, a 4-dimensional complex vector space.

Recall also that the four sphere $S^4$ is naturally embedded into $\m$.
It can be realized as the fixed point set of the following real structure
$\sigma$ on $\m$:
\begin{equation} \label{m.invol}
\begin{array}{ccc}
\sigma(z_{11'}) = \overline{z_{22'}} & \ \ \ & \sigma(z_{12'}) = -\overline{z_{21'}} \\
\sigma(z_{21'}) = -\overline{z_{12'}}& \ \ \ & \sigma(z_{22'}) = \overline{z_{11'}}\\
\sigma(D) = \overline{D} & \ \ \ & \sigma(D') = \overline{D'}
\end{array} \end{equation}
Moreover, we observe that even though the affine pieces $\m^{\rm I}$ and 
$\m^{\rm J}$ do not cover $\m$, all of its real points lie within their
union, that is $S^4\hookrightarrow\m^{\rm I}\cup\m^{\rm J}$.

We claim that there are $q$-deformations of
$\mathbf{G}_2(\st)$ and $\mathbf{F}_{1,2}(\st)$
which provide a correspondence between a $q$-deformed
Grassmannian and the {\em classical} twistor space $\mathbf{P}(\st)$. 
Furthermore, our quantum Minkowski space-time $\cm^{\rm I}_q$
is an {\em affine patch} of the $q$-deformed Grassmannian.
Moreover, all relevant noncommutative varieties can also be
obtained from the quantum group $GL(4)_q$ extended by appropriate
derivations. However, in the construction below, we will only use
the quantum group $SL(2)_q$ enlarged by its corepresentations of
functional dimension 2 and by certain degree operators.

Once these noncommutative varieties are constructed, we hope that
our conjecture will be proved in complete parallel with 
the classical version.

\paragraph{Quantum Grassmannian}
Let us now describe the noncommutative variety from which the
quantum Minkowski space-time $\cm^{\rm I}_q$ is obtained via localization.
Recall that $q$ is a formal parameter.

\begin{definition} \label{mpq}
The quantum compactified, complexified Minkowski space $\cm_{p,q}$ 
is the associative graded $\cpx$-algebra generated by 
$z_{11'},z_{12'},z_{21'},z_{22'},D,D'$ satisfying the relations
(\ref{relations1}) to (\ref{relations5}) below $(p=q^{\pm 1})$:
\end{definition}

\begin{equation} \label{relations1}
\begin{array}{lcr}
z_{11'}z_{12'}=z_{12'}z_{11'} & \ \ \ & z_{11'}z_{21'}=z_{21'}z_{11'} \\
z_{12'}z_{22'}=z_{22'}z_{12'} & \ \ \ & z_{21'}z_{22'}=z_{22'}z_{21'} \\
& z_{12'}z_{21'}=z_{21'}z_{12'} &
\end{array} \end{equation}
\begin{equation} \label{relations2}
q^{-1} (z_{11'}z_{22'}-z_{12'}z_{21'})=q (z_{22'}z_{11'}-z_{12'}z_{21'})
\end{equation}
\begin{equation} \label{relations3}
\begin{array}{lcr}
Dz_{11'}=pq^{-1} z_{11'}D & \ \ \ & D'z_{11'}=p^{-1}q^{-1} z_{11'}D' \\
Dz_{12'}=pq^{-1} z_{12'}D & \ \ \ & D'z_{12'}=p^{-1}q z_{12'}D'      \\
Dz_{21'}=pq z_{21'}D      & \ \ \ & D'z_{21'}=p^{-1}q^{-1} z_{21'}D' \\
Dz_{22'}=pq z_{22'}D      & \ \ \ & D'z_{22'}=p^{-1}q z_{22'}D'      \\
\end{array} \end{equation}
\begin{equation} \label{relations4}
p^{-1} DD'=p D'D
\end{equation}
\begin{equation} \label{relations5}
q^{-1} (z_{11'}z_{22'}-z_{12'}z_{21'}) = p^{-1} DD'
\end{equation}

The relations (\ref{relations1})-(\ref{relations4}) are simply
commutation relations, while (\ref{relations5}) plays the role
of the quadric (\ref{quad}) that defines $\m$ as a subvariety
of $\p^5$. In other words, the algebra $\cm_{p,q}$ can be regarded
as a {\em quantum Grassmannian}. Note also that the relations
(\ref{relations1})-(\ref{relations5}) can be expressed in $R$-matrix form.

Furthermore, in analogy with the classical case, it is not difficult
to establish the following:
 
\begin{proposition} \label{local}
The algebras $\cm^{\rm I}_q$ and $\cm^{\rm J}_q$ are localizations of
$\cm_{p,q}$ with respect to $D$ and $D'$, respectively, independently
on $p=q^{\pm1}$. In particular,
$$ x_{rs'} = \frac{z_{rs'}}{D} \ \ \ {\rm and} \ \ \ y_{rs'} =
\frac{z_{rs'}}{D'} $$
\end{proposition}
\begin{proof}
The proof is a straightforward calculation left to the reader. Note only
that the notation $\frac{z_{rs'}}{D}$ means 
$c_{rs'}^{-1/2}D^{-1}z_{rs'} = c_{rs'}^{1/2}z_{rs'}D^{-1}$ whenever
$D^{-1}z_{rs'}=c_{rs'}z_{ij'}D^{-1}$ for some constant $c_{rs'}$, and
similarly for $\frac{z_{rs'}}{D'}$.
\hfill $\Box$ \end{proof}

It is also important to note that it follows from the comparison
of equations (\ref{relations3}) with equations 
(\ref{xvar}) and (\ref{yvar}) that $\delta^2=\frac{D'}{D}$.

In geometric terms, our quantum Minkowski space-time $\cm^{\rm I}_q$
appears as an affine patch of the quantum Grassmannian $\cm_{p,q}$, 
thus justifying our Definition.

\paragraph{The quantum 4-sphere}
Let us define the action of the $\dagger$-involution on the
generators of $\cm_{p,q}$ as follows:
$$ \begin{array}{lcl}
z_{11'}^\dagger = z_{22'} & \ \ \ & z_{12'}^\dagger = - z_{21'} \\
z_{21'}^\dagger = - z_{12'} & \ \ \ & z_{22'}^\dagger = z_{11'} \\
D^\dagger = D & \ \ \ & D'^\dagger = D'
\\
q^\dagger = q & \ \ \ & p^\dagger = p
\end{array} $$
We then extend $\dagger$ to 
$\cm_{p,q}$ by requiring it
to be  a conjugate linear anti-homomor-phism.
Note that $\dagger$ interchanges the quantum Grassmannians, i.e. 
$\dagger:\cm_{p,q}\to\cm_{p^{-1},q}$. One easily checks that
it is consistent with the $\dagger$-involutions
previously defined on
$\cm^{\rm I}_q$ and $\cm^{\rm J}_q$.
Also, comparing with
(\ref{m.invol}), we see that the map $\dagger$
just defined is the analogue of
the real structure $\sigma$ on $\m$.
Therefore, we propose the following definition. 

\begin{definition} \label{spq}
The quantum 4-dimensional sphere $S_{q}$ is the pair of algebras 
$\cm_{q^{\pm1},q}$ equipped with the map $\dagger:\cm_{p,q}\to\cm_{p^{-1},q}$: 
$$ S_q = (\cm_{q^{\pm1},q},\dagger) $$
\end{definition}

Alternative definitions of quantum 4-dimensional spheres have been
proposed by several other authors, see in particular \cite{BCT,CL,DLM}.
Note however that $\cm_{p,q}$ is a deformation of the algebra of 
homogeneous coordinates, not of the algebra of polynomial functions
as in the examples discussed in the references mentioned above.
Our definition is justified by the analogy with the commutative case
and by the construction of quantum instantons
on the two ``affine patches'' $S^{\rm I}_q$ and $S^{\rm J}_q$ of
$S_q$ done in Section 3.

\paragraph{Quantum flag variety}
Let $\p$ denote the projective twistor space $\mathbf{P}(\st)$;
recall that the flag variety $\f=\mathbf{F}_{1,2}(\st)$ can be naturally
embedded in the product $\p\times\m$ by sending the flag 
[line$\subset$plane]$\in\f$ to the pair
(line,plane)$\in\p\times\m$.  More precisely, let $[z_k]$
($k=1,2,1',2'$, consistently with the splitting (\ref{decomp}))
be homogeneous coordinates in $\p$; then $\f$ can be described
as an intersection of the following quadrics in $\p\times\m$:
\begin{eqnarray}
Dz_{1'} + z_1z_{21'} - z_2z_{11'} = 0 & \ \ \ &  
Dz_{2'} + z_1z_{22'} - z_2z_{12'} = 0\label{q1} \\
D'z_1 + z_{2'}z_{11'} - z_{1'}z_{12'} = 0 & \ \ \ &
D'z_2 + z_{2'}z_{21'} - z_{1'}z_{22'} = 0 \label{q2}
\end{eqnarray}

Let us now introduce the quantum flag variety that provides a correspondence
between the quantum Grassmannian described above and the classical
twistor space. Notice that the coordinate algebra of the projective 
twistor space $\p$ is simply given by:
$$ \cp  = \cpx[z_1,z_2,z_{1'},z_{2'}]_h $$
where the subscript ``$h$" means that $\cp$ consists only of 
the homogeneous polynomials. 

Again, the starting point lies in a natural extension of the quantum
group $SU(2)_q$, this time by left and right corepresentations of
functional dimension 2. We define left and right corepresentations of
$GL(2)_q$ as the vector spaces $L_q$ and $R_q$ of polynomials on two
noncommutative variables $g_1,g_2$ and $g_{1'},g_{2'}$, respectively,
satisfying the commutation relations:
\begin{equation} \label{left}
g_1g_2 = q^{-1}g_2g_1
\end{equation}
\begin{equation} \label{right}
g_{1'}g_{2'} = q^{-1}g_{2'}g_{1'}
\end{equation}
The coaction is given by, respectively:
\begin{equation} \label{lcoation}
\Delta_L(g_r) = \sum_{k=1,2} g_{rk'}\otimes g_k
\end{equation}
\begin{equation} \label{rcoation}
\Delta_R(g_{s'}) = \sum_{k=1,2} g_{k'}\otimes g_{ks'}
\end{equation}
The commutation relations (\ref{gl2q1}) and (\ref{gl2q2}) for the
quantum group $GL(2)_q$ imply that $\Delta_L$ and $\Delta_R$ are 
indeed corepresentations (in fact, they are also necessary conditions).

We can also regard $L_q$ and $R_q$ as corepresentations of the quantum
groups $SL(2)_q$, $\widetilde{GL}(2)_q$ and $\widetilde{SL}(2)_q$.
Again, we will be primarily interested in $\widetilde{SL}(2)_q$.

We define the $\dagger$-involution on the generators of $L_q$ and $R_q$ as
follows:
\begin{equation} \label{co.star}
\begin{array} {ccc}
g_1^\dagger  = g_2, & & g_2^\dagger = -g_1 \\
g_{1'}^\dagger = g_{2'}, & & g_{2'}^\dagger = -g_{1'}
\end{array} \end{equation}
As before, we extend $\dagger$ to a conjugate linear anti-homomorphisms of
$L_q$ and
$R_q$; see page \pageref{ext.star}. We also get:
\begin{equation} \label{corep.star}
\Delta_L(x^\dagger) = (\Delta_L(x))^\dagger\otimes(\Delta_L(x))^\dagger,
\ \ \ {\rm and}\ \ \ 
\Delta_R(y^\dagger) = (\Delta_R(y))^\dagger\otimes(\Delta_R(y))^\dagger
\end{equation}
for $x\in L_q$ and $y\in R_q$, which means that $\Delta_L$ and 
$\Delta_R$ are well defined as corepresentations of the quantum 
groups $SU(2)_q$ and $\widetilde{SU}(2)_q$.

Next, we define the semi-direct product quantum groups 
$(SL(2)_q\ltimes L_q)_{p}$ and $(R_q\rtimes SL(2)_q)_{p}$ as the algebras
with the combined generators of $SL(2)_q$ and $L_q$ and $SL(2)_q$ 
and $R_q$, respectively, satisfying the relations (in $R$-matrix form):
$$ \begin{array}{ccc}
p^{1/2} R T_1 S_2 = S_2 T_1 &
\ \ \ {\rm and}\ \ \ &
T_1 S_2^\prime = S_2^\prime T_1 R p^{1/2}
\end{array} $$
where $R$ is the $R$-matrix given in (\ref{rmatrix}) and
$$ \begin{array}{ccc}
S = \left( \begin{array}{c} g_1 \\ g_2 \end{array} \right) & 
\ \ \ {\rm and}\ \ \ &
S' = \left( \begin{array}{c} g_{1'} \\ g_{2'} \end{array} \right)
\end{array} $$
are the generating matrices for $L_q$ and $R_q$, respectively, 
with $S_2=\id\otimes S$, $S_2'=\id\otimes S'$. We define
a comultiplication on $(SL(2)_q\ltimes L_q)_{p}$ and $(R_q\rtimes SL(2)_q)_{p}$ 
by the formulas (\ref{comult}), (\ref{lcoation}) and
(\ref{rcoation}). Remark that:

\begin{proposition} \label{cons}
The multiplication and comultiplication in $(SL(2)_q\ltimes L_q)_{p}$ and 
$(R_q\rtimes SL(2)_q)_{p}$ are consistent, and they yield a bialgebra 
structure.
\end{proposition}

The bialgebras
$(SL(2)_q\ltimes L_q)_{p}$ and $(R_q\rtimes SL(2)_q)_{p}$ are quantum
versions of the classical groups $SL(2)\ltimes L$ and $R\rtimes SL(2)$,
respectively, where $L$ and $R$ denote left and right two-dimensional
representations of $SL(2)$. It is well known that these two groups are
isomorphic, where the isomorphism is given by the identification
$L\simeq R$. It turns out that this isomorphism has a quantum analogue,
and the bialgebras $(SL(2)_q\ltimes L_q)_{p}$ and  $(R_q\rtimes SL(2)_q)_{p}$
are also isomorphic. The quantum isomorphism 
$\kappa:(SL(2)_q\ltimes L_q)_{p}\seta(R_q\rtimes SL(2)_q)_{p}$  is the
identity on $SL(2)_q$ and acts in the generators of $L_q$ as follows:
\begin{eqnarray*}
\kappa(g_1) & = & 
p^{-3/4} \left( -q^{1/2}g_{11'}g_{2'} + q^{-1/2}g_{12'}g_{1'} \right) \\
\kappa(g_2) & = & 
p^{-3/4} \left( -q^{1/2}g_{21'}g_{2'} + q^{-1/2}g_{22'}g_{1'} \right)
\end{eqnarray*}
Therefore we define $F_{p,q}$ as $R_q\rtimes SL(2)_q\ltimes L_q$ modulo
the equivalence given by the isomorphism $\kappa$. One can compute the
following commutation relations:
\begin{equation} \label{g.commut} \begin{array}{ccc}
g_1g_{1'} = p^{1/2} g_{1'}g_1 & \ \ \ & g_2g_{1'} = p^{1/2} g_{1'}g_2 \\
g_1g_{2'} = p^{1/2} g_{2'}g_1 & \ \ \ & g_2g_{2'} = p^{1/2} g_{2'}g_2
\end{array} \end{equation}

We proceed by adjoining two new generators $\partial$ and $\partial'$
to $F_{p,q}$, and postulating the following relations
(keeping in mind that $p=q^{\pm1}$):
\begin{equation} \label{dd'grs}  \begin{array}{lcl} 
\partial g_{11'} = q^{-1/2} g_{11'} \partial & \ \ \ & 
\partial' g_{11'} = q^{-1/2}  g_{11'} \partial' \\
\partial g_{12'} = q^{-1/2} g_{12'} \partial & \ \ \ & 
\partial' g_{12'} = q^{1/2} g_{12'} \partial' \\
\partial g_{21'} = q^{1/2} g_{21'} \partial & \ \ \ & 
\partial' g_{21'} = q^{-1/2}  g_{21'} \partial' \\
\partial g_{22'} = q^{1/2} g_{22'} \partial & \ \ \ & 
\partial' g_{22'} = q^{1/2} g_{22'} \partial' \\
\end{array} \end{equation}
\begin{equation} \label{dd'gk}  \begin{array}{lcl} 
\partial g_1 = p^{1/2}q^{-1/2} g_1 \partial & \ \ \ & 
\partial' g_1 = p^{1/2}  g_1 \partial' \\
\partial g_2 = p^{-1/2}q^{1/2} g_2 \partial & \ \ \ & 
\partial' g_2 = p^{1/2} g_2 \partial' \\
\partial g_{1'} = p^{-1/2} g_{1'} \partial & \ \ \ & 
\partial' g_{1'} = p^{1/2}q^{-1/2}  g_{1'} \partial' \\
\partial g_{2'} = p^{-1/2} g_{2'} \partial & \ \ \ & 
\partial' g_{2'} = p^{1/2}q^{1/2} g_{2'} \partial' \\
\end{array} \end{equation}
\begin{equation} \label{dd'}
\partial \partial' =  p^{1/2} \partial' \partial
\end{equation}
We denote the extended algebra $F_{p,q}[\partial,\partial']$ by
$\widetilde{F}_{p,q}$; this is the prototype of our quantum flag
variety.

Furthermore, we also define the $\dagger$-involution on $\widetilde{F}_{p,q}$ 
by extending the $\dagger$-involutions on $\widetilde{SL}(2)_q$, $L_q$ and $R_q$
defined on (\ref{involution}) and (\ref{co.star}) and declaring:
\begin{equation} \label{dd'star}
\partial^\dagger =\partial 
\ \ \ {\rm and} \ \ \ 
\partial'^\dagger =  \partial'
\end{equation}
 
As in the construction of the quantum Minkowski space-time,
we introduce a new set of variables $z_{rs'},z_k,z_{k'}$ 
in the bialgebra $\widetilde{F}_{p,q}=F_{p,q}[\partial,\partial']$,
related to the original variables $g_{rs'},g_k,g_{k'}$ via:
$$ \frac{z_{rs'}}{\Delta} = g_{rs'} \ \ \ \ \ 
\frac{z_k}{\partial} = g_k \ \ \ \ \ 
\frac{z_{k'}}{\partial'} = g_{k'} $$
where $\Delta=p^{1/4} \partial' \partial=p^{-1/4} \partial \partial'$.
The fractions above have the same meaning as in Proposition \ref{local}.
We also set: 
$$ D = \partial^2 \ \ \ \ \ D' = \partial'^2 $$

Using relations (\ref{dd'grs}) and (\ref{dd'}), one re-obtains the
relations (\ref{relations1}-\ref{relations5}), which define
$\cm_{p,q}$. In other words, the quantum compactified, complexified 
Minkowski space $\cm_{p,q}$ can be regarded as a subalgebra of
the extended quantum group $SL(2)_q[\partial,\partial']$. Notice
in particular that the quantum quadric (\ref{relations5}) is just
the relation $\det_q(T)=1$ expressed in the new variables $z_{rs'}$.

Now using relations (\ref{dd'gk}) and (\ref{dd'}), it is easy to check that:
\begin{equation} \label{relations6.5}
z_1,z_2,z_{1'}z_{2'}\ \text{commute with one another}
\end{equation}
Furthermore, we also obtain:
\begin{equation} \label{relations7}
\begin{array}{lcr}
z_1z_{11'}=z_{11'}z_1       & \ \ \ & z_1z_{12'}=z_{12'}z_1      \\
z_2z_{21'}=z_{21'}z_2       & \ \ \ & z_2z_{22'}=z_{22'}z_2      \\
z_{1'}z_{11'}=z_{11'}z_{1'} & \ \ \ & z_{1'}z_{21'}=z_{21'}z_{1'} \\
z_{2'}z_{12'}=z_{12'}z_{2'} & \ \ \ & z_{2'}z_{22'}=z_{22'}z_{2'}
\end{array} \end{equation}
and, keeping in mind that $p=q^{\pm1}$:
\begin{equation} \label{relations8}
\begin{array}{c}
z_1z_{21'} = pq z_{21'}z_1 + (1-pq) z_{11'}z_2 \\
z_1z_{22'} = pq  z_{22'}z_1 + (1-pq)  z_{12'}z_2 \\
z_2z_{11'} = pq^{-1}  z_{11'}z_{2} + (1-pq^{-1})  z_{21'}z_1 \\
z_2z_{12'} = pq^{-1}  z_{22'}z_{2} + (1-pq^{-1})  z_{22'}z_1 \\
z_{1'}z_{12'} = p^{-1}q  z_{12'}z_{1'} + (1-p^{-1}q)  z_{11'}z_{2'} \\
z_{1'}z_{12'} = p^{-1}q  z_{22'}z_{1'}+ (1-p^{-1}q)  z_{21'}z_{2'} \\
z_{2'}z_{11'} = p^{-1}q^{-1}  z_{11'}z_{2'} + (1-p^{-1}q^{-1})  z_{12'}z_{1'} \\
z_{2'}z_{21'} = p^{-1}q^{-1}  z_{21'}z_{2'}+ (1-p^{-1}q^{-1})  z_{22'}z_{1'}
\end{array} \end{equation}
The commutation relations between the degree operators $D$ and $D'$
and the commuting generators $z_1,z_2,z_{1'},z_{2'}$ are given by:
\begin{equation} \label{relations9}
\begin{array}{lcl}
Dz_1= p^{-1}q^{-1}  z_1D & \ \ \ & D'z_1=z_1D' \\
Dz_2= p^{-1}q  z_2D      & \ \ \ & D'z_2=z_2D' \\
Dz_{1'}=z_{1'}D              & \ \ \ & D'z_{1'}=pq^{-1}  z_{1'}D'\\
Dz_{2'}=z_{2'}D              & \ \ \ & D'z_{2'}=pq  z_{2'}D'
\end{array} \end{equation}
The last set of relations below follows from the identifications induced
by the isomorphism $\kappa$ and its inverse, and yields the quantum
analogues of the quadrics (\ref{q1}) and (\ref{q2}) defining $\f$ as a 
subvariety of $\p\times\m$:
\begin{equation} \label{relations10}
\begin{array}{c}
Dz_{1'} = p      ( z_{11'}z_2    - z_{21'}z_{1} ) \\
Dz_{2'} = p      ( z_{12'}z_2    - z_{22'}z_{1} ) \\
D'z_1   = p^{-1} (-z_{11'}z_{2'} + z_{12'}z_{1'}) \\
D'z_2   = p^{-1} (-z_{21'}z_{2'} + z_{22'}z_{1'})
\end{array} \end{equation}
This observation motivates our next definition:

\begin{definition} \label{qflag}
The quantum flag variety $\cf_{p,q}$ is the associative graded $\cpx$-algebra
with generators $z_{11'},z_{12'},z_{21'},z_{22'},D,D',z_1,z_2,z_{1'}z_{2'}$ 
satisfying relations (\ref{relations1})-(\ref{relations5}) and
(\ref{relations7})-(\ref{relations10}) above.
\end{definition}

In particular, note that $\cf_{p,q}$ can be regarded as a 
subalgebra of $\widetilde{F}_{p,q}$.
We also remark that:
\begin{equation} \label{qflag2}
\cf_{p,q}=\cm_{p,q}\otimes_{\cpx}\cp/(\ref{relations7})-(\ref{relations10})
\end{equation}
in close analogy with the classical case. All the relations
(\ref{relations7}-\ref{relations10}) can also
be expressed in $R$-matrix form.

Now let 
$$ \cf_q^{\rm I} = \cf_{p,q}[D^{-1}]_0 \simeq 
\cm_q^{\rm I}\otimes\cpx[z_1,z_2]_h $$ 
$$\cf_q^{\rm J}=\cf_{p,q}[D'^{-1}]_0 \simeq 
\cm_q^{\rm J}\otimes\cpx[z_1,z_2]_h $$ 
where $\cf_{p,q}[D^{-1}]$ and $\cf_{p,q}[D'^{-1}]$ denote the localization of
$\cf_{p,q}$ as a $\cm_{p,q}$-bimodule and 
the subscript ``$0$" means that we take only the degree zero part of
the localized graded algebra, and the subscript ``$h$" means that 
only the homogeneous polynomials are considered. Geometrically, notice
that these
algebras are playing the roles of $q$-deformations of the ``affine''
flag varieties $\f^{\rm I}=\nu^{-1}(\m^{\rm I})=\m^{\rm I}\times\p^1$ and 
$\f^{\rm J}=\nu^{-1}(\m^{\rm J})=\m^{\rm J}\times\p^1$
\cite{WW}. To further justify Definition \ref{qflag} we prove:

\begin{theorem} \label{nc.diagram}
The maps $\sm_{p,q}:\cp\seta\cf_{p,q}$ and
$\sn_{p,q}:\cm_{p,q}\seta\cf_{p,q}$ defined as identities on the
generators are injective homomorphisms. Furthermore, the 
multiplication in the associative algebras $\cm_{p,q}$, 
$\cm_{p,q}^{\rm I}$,  $\cm_{p,q}^{\rm J}$ and $\cf_{p,q}$,
$\cf_{p,q}^{\rm I}$,  $\cf_{p,q}^{\rm J}$  
described above is consistent, and they have the same Hilbert
polynomials as their commutative counterparts $\m$, $\m^{\rm I}$,
$\m^{\rm J}$ and $\f$, $\f^{\rm I}$,  $\f^{\rm J}$, respectively.
\end{theorem}

\begin{proof} 
The statement regarding the maps $\sm_{p,q}$ and
$\sn_{p,q}$ is clear from our construction, see (\ref{qflag2}). 

For the second part of the Theorem, we must first check the 
consistency of multiplication in the algebra $\cm_{p,q}$ by
embedding it into its localization $\cm_{p,q}[D^{-1}]$. The commutation
relations (\ref{relations1}), (\ref{relations2}) and the first column
of (\ref{relations3}) allow us to choose a basis in 
$\cpx[z_{rs'},D^{\pm1}]$ of the form:
\begin{equation} \label{basis}
z_{11'}^{n_{11'}}z_{12'}^{n_{12'}}z_{21'}^{n_{21'}}z_{22'}^{n_{22'}}D^n,
\ \ \ {\rm with} \ \ \ n_{rs'}\in\mathbb{Z}_+,\ \ n\in\mathbb{Z}
\end{equation}
Inverting $D$ in (\ref{relations5}) we obtain an expression for $D'$,
and we can check directly the second column of (\ref{relations3})
and relation (\ref{relations4}). Thus the basis (\ref{basis}) is in
fact a basis of $\cm_{p,q}[D^{-1}]$, and its Hilbert polynomial
coincides with the classical one. 

Next we find the Hilbert polynomial of the algebra $\cm_{p,q}$.
Using (\ref{relations5}), we can present an arbitrary element from
$\cm_{p,q}$ uniquely in the following form:
\begin{equation} \label{poly}
P_0(z_{rs'}) + P_1(z_{rs'},D)D + P_2(z_{rs'},D')D'
\end{equation}
where $P_0,P_1,P_2$ are polynomials and $z_{rs'}$ are ordered as in 
(\ref{basis}). This presentation coincides with the classical one,
and allows to compute explicitly the Hilbert polynomial of
$\cm_{p,q}$.

Similarly, we check the consistency of multiplication in the algebra
$\cf_{p,q}$ by embedding it into $\cf_{p,q}[D^{-1}]$. In addition to the
relations used above, we use the the first half of the commutation relations
(\ref{relations7}), (\ref{relations8}) and the first column of
(\ref{relations9}) to choose a basis in $\cpx[z_{rs'},z_{k},D^{\pm1}]$
of the form:
\begin{equation} \label{basis2}
z_{11'}^{n_{11'}}z_{12'}^{n_{12'}}z_{21'}^{n_{21'}}z_{22'}^{n_{22'}}
z_1^{n_1}z_2^{n_2}D^n,
\ \ \ {\rm with} \ \ \ n_{rs'},n_k\in\mathbb{Z}_+,\ \ n\in\mathbb{Z}
\end{equation}
Inverting $D$ in (\ref{relations5}) and also in the first half of
(\ref{relations10}), we obtain expressions for $D'$, $z_{1'}$ and
$z_{2'}$ and check directly all the other relations involving these
three generators, namely the second half of the commutation relations
(\ref{relations7}), (\ref{relations8}) and the second column of
(\ref{relations9}). We conclude that the multiplication in
$\cf_{p,q}[D^{-1}]$ is consistent and its basis is given by 
(\ref{basis2}).

Finally, we obtain the Hilbert polynomial of $\cf_{p,q}$ by
noticing that using (\ref{relations5}) and (\ref{relations10})
we can present an arbitrary element of $\cf_{p,q}$ uniquely
in the form
\begin{equation} \label{poly2}
P_0(z_{rs'},z_k,z_{k'}) + P_1(z_{rs'},z_k,D)D + P_2(z_{rs'},z_{k'},D')D'
\end{equation}
where $P_0,P_1,P_2$ are polynomials and $z_{rs'},z_k,z_{k'}$ are ordered
as in (\ref{basis2}). Again this presentation coincides with the
classical one, and yields an explicit formula for the Hilbert
polynomial of $\cf_{p,q}$.

Alternatively, one could also check the consistency of multiplication
more efficiently using the $R$-matrix formulation.
\hfill $\Box$ \end{proof}

Summing up, we have constructed noncommutative varieties $\cm_{p,q}$
and $\cf_{p,q}$, thought as quantum deformations of $\m$ and $\f$, 
and injective morphisms $\sm_{p,q}$ and $\sn_{p,q}$ fitting into the 
following diagram:
\begin{equation} \label{nc-dblfib}
\xymatrix{ & \cf_{p,q} & \\ 
\cp \ar[ur]_{\sm_{p,q}} & & \cm_{p,q} \ar[ul]^{\sn_{p,q}} }
\end{equation}
where $\cp$ is just the commutative projective 3-space. In analogy with
the classical case, diagram (\ref{nc-dblfib}) will be called the 
{\em quantum twistor diagram}. 

Furthermore, we also have the {\em quantum local twistor
diagrams}:
\begin{equation} \label{loc.nc-dblfib} \begin{array}{lr}
\xymatrix{ & \cf_q^{\rm I} & \\ 
\cp^{\rm I} \ar[ur] & & \cm_q^{\rm I} \ar[ul] } &
\xymatrix{ & \cf_q^{\rm J} & \\ 
\cp^{\rm J} \ar[ur] & & \cm_q^{\rm J} \ar[ul] }
\end{array} \end{equation}
where 
$$ \cp^{\rm I} = \cpx[z_1 , z_2 , 
x_{11'}z_2-x_{21'}z_1 , x_{12'}z_2-x_{22'}z_1]_h $$
$$ \cp^{\rm J} = \cpx[ y_{12'}z_{1'}-y_{11'}z_{2'} , 
y_{22'}z_{1'}-y_{21'}z_{2'} , z_{1'} , z_{2'}]_h $$
as commutative subalgebras of $\cf^{\rm I}_q$ and $\cf^{\rm J}_q$, 
respectively. Notice that $\cp^{\rm I}$ and $\cp^{\rm J}$
are exactly the coordinate algebras of
$\mu(\nu^{-1}(\m^{\rm I}))$ and $\mu(\nu^{-1}(\m^{\rm J}))$, the
``affine'' twistor spaces \cite{WW}.

It is also clear from our construction that the quantum 
twistor diagram (\ref{nc-dblfib}) comes with a {\em real structure}
defined by the $\dagger$ involutions on $\cf_{p,q}$ and $\cm_{p,q}$
defined above, acting on the generators of $\cp$ as follows:
$$ z_1^\dagger = z_2  \ \ \ z_2^\dagger = - z_1
\ \ \  z_{1'}^\dagger = z_{2'} \ \ \ z_{2'}^\dagger = - z_{1'} $$

After rewriting the classical Penrose transform in algebraic terms, 
we can use the quantum twistor diagrams (\ref{nc-dblfib})
and (\ref{loc.nc-dblfib})
to transform {\em quantum} objects (solutions of the quantum
ASDYM equation) into {\em classical} ones (holomorphic bundles with a
real structure). This is the strategy for the proof of our completeness
conjecture, and a topic for a future paper. 

%----------------------------------

\subsection{Remarks on roots of unity and representation theory}
\label{ru}

In our construction of quantum instantons, we have specialized the
formal
parameter $q$ to a positive real number. Thus it is natural
to ask what happens if we consider other specializations. Since our
construction is deeply related with the behaviour of the Haar functional
on the quantum group $SL(2)_q$, we expect that our results remain valid
for all specializations of $q$ for
which the Haar functional has a similar structure,
namely for any complex $q$ that is not a root of unity.

It is an interesting problem
to obtain a modified ADHM data for any
given root of unity, and
then describe the corresponding moduli space
of quantum instantons in
terms of this data. We expect that such
description will be given by
certain quiver varieties.

Moreover, equivariant anti-self-dual connections on $\real^4$ has long been
a topic of intense research. In particular, instantons on the quotient
spaces $\cpx^2/\Gamma$, where $\Gamma$ is a finite subgroup of $SU(2)$,
and on their desingularizations (known as ALE spaces) have been studied
by Kronheimer and Nakajima in \cite{KN,N4}. The moduli spaces of instantons
on ALE spaces play a fundamental role on Nakajima's construction of
representations of affine Lie algebras \cite{N3}.

To obtain a quantum analogue of this construction, we first need to
describe ``finite subgroups'' of $SU(2)_q$. This turns out to be a deep
problem, recently solved by Ostrik \cite{O}, though only at the
level of representation categories. The resulting classification, as
in the classical case, is given by Dynkin diagrams of ADE type. This
suggests that the corresponding moduli spaces of quantum instantons
might shed a new light at Nakajima's quiver varieties and his results
on representation theory. 

We would also like to notice the similarities between the deformation of
the representation theory for $SU(2)$ (as well as other compact simple
Lie groups) and the deformation of instanton theory presented in this
paper.
As it is well known, the classification of irreducible 
corepresentations
of the quantum group $SU(2)_q$, the behaviour of the
Haar functional, the decomposition of their
tensor products, etc, 
is the same as in the classical, non-deformed case if
$q$ is not a root of unity, see \cite{CP}. 

The anti-self-duality condition can be viewed as {\em half} of the relations
necessary to determine a representation of the quantum group $SU(2)_q$, and
we also obtain a classification identical to the non-deformed case.
Tensor products of representations of quantum groups correspond in our
context to the composition of instantons of different ranks, and again
we expect it not to depend on the quantization parameter, if it is not
a root of unity.

%----------------------------------

\subsection{Further perspectives}

It is well known that the Penrose transform can be used to obtain
solutions not only of the ASDYM equation, but also of a variety of
other differential equations on Minkowski space-time \cite{WW}. One can then
expect to define the quantum analogue of such equations in a way that our
proposed quantum Penrose transform still yields the same result
for generic $q$.

The most basic class of examples are the massless free field equations,
which includes the wave, Dirac and Maxwell equations. In particular, the
polynomial solutions to the wave equation on $\mi$ are given by the
matrix coefficients of irreducible representations of $SU(2)$. Thus the
appropriate quantum wave equation should have solutions on $\cm^{\rm I}_q$
given by the matrix coefficients of representations of the quantum group
$SU(2)_q$, and a similar relation should also be true for the higher spin
equations.

Moreover, the wave equation for the anti-self-dual connection plays a
prominent role in the proof of completeness of instantons \cite{Ma} and
its solutions correspond, via Penrose transform, to certain sheaf
cohomology classes defined over $\mathbf{P}(\st)$. We expect that all this
structure will be preserved under our quantum deformation of the
twistor diagram, leading to a proof of our completeness conjecture.

Finally, we cannot avoid to mention another remarkable noncommutative
deformation of instantons discovered by Nekrasov and Schwarz \cite{NS}.
The moduli spaces of such noncommutative instantons can
 be regarded as
a smooth compactification of the classical ADHM data, though the
relation with the latter is not as straightforward as in our case \cite{N1}.

A version of the Penrose transform for the Nekrasov-Schwarz
noncommutative instantons has been discussed in \cite{KKO}. Also, it
has been observed in \cite{L} that the factorization of the
noncommutative Minkowski space by a finite subgroup
$\Gamma\hookrightarrow SU(2)$ in the Nekrasov-Schwarz setting yields the
smooth varieties previously introduced by Nakajima \cite{N3}.

It is an interesting and challenging problem to obtain a natural
compactification of the moduli spaces of quantum instantons by
considering a more general type of $\cm^{\rm I}_q$-modules than those
considered here, and then compare with the smooth compactification
appearing in the Nekrasov-Schwarz approach.
This might lead to a further extension of our understanding of the
mathematical structure of space-time.

%---------------------------------------------------------------
%---------------------------------------------------------------

\paragraph{Acknowledgement.}
We thank Yu. Berest for a discussion. We also thank
the referee for his comments on the first version
of this paper.
 
I.F. is supported by the NSF grant DMS-0070551, and
would like to thank the MSRI for its hospitality at
the later stages of this project.
M.J. thanks the Departments of Mathematics at
Yale University and at the University of Pennsylvania.

 \end{document}